\documentclass[a4paper,twoside,10pt]{article}

\usepackage[a4paper,width=150mm,top=25mm,bottom=25mm,bindingoffset=6mm]{geometry}

\usepackage[english]{babel}

\usepackage{graphicx}
\usepackage[font=small,labelfont=bf]{caption}
\usepackage{subcaption}
\usepackage{color}
\usepackage{mdframed}
\usepackage{framed} 
\usepackage{cancel}
\usepackage[normalem]{ulem}

\usepackage{lineno,hyperref}
\modulolinenumbers[5]

\usepackage{amsmath,amssymb,amsthm,mathtools,mathrsfs,epsfig,amsfonts,wasysym}
\usepackage{MnSymbol}

\usepackage[normalem]{ulem}
\usepackage{ragged2e}
\usepackage{csquotes}
\usepackage{tasks}
\usepackage{paralist}
\usepackage[inline]{enumitem}
\usepackage{tikz-cd}
\usepackage{stmaryrd}

\usepackage[titletoc,title]{appendix}

\newcommand{\de}{\mathrm{d}}

\usepackage{url,hyperref}
\hypersetup{colorlinks=true,linkcolor=blue,filecolor=magenta,urlcolor=cyan}
\usepackage{nicematrix}
\usepackage{tikz}
\usetikzlibrary{fit}

\usepackage[all]{xy}

\newtheorem*{example*}{Example}
\newtheorem*{reminder*}{Reminder}
\newtheorem*{remark*}{Remark}
\newtheorem*{note*}{Note}

\usepackage{authblk}
\usepackage{todonotes}

\definecolor{wildstrawberry}{rgb}{1.0, 0.26, 0.64}

\definecolor{ao(english)}{rgb}{0.0, 0.5, 0.0}

\begin{document}

\title{On the quasi-steady-state assumption in enzyme kinetics: rigorous analysis\vspace{.5cm}\\}

\author[1]{Vasiliki Bitsouni\thanks{\texttt{vbitsouni@math.uoa.gr}}}
\author[1]{Nikolaos Gialelis\thanks{\texttt{ngialelis@math.uoa.gr}}}
\author[1]{Ioannis G. Stratis\thanks{\texttt{istratis@math.uoa.gr}}}
\affil[1]{Department of Mathematics,\\ National and Kapodistrian University of Athens, Panepistimioupolis, GR-15784 Athens, Greece}

\date{}

\maketitle

\begin{abstract}
We study, from a purely quantitative point of view, the quasi-steady-state assumption for the fundamental mathematical model of the general enzymatic reaction: we re-establish, on a rigorous basis, certain already known results and we propose a novel approach to the subject, that makes waivable an ambiguous, widely used, practice. In particular, we define the two parts of the assumption in a quantitative fashion, we employ a simple algorithm for the proper scaling of the corresponding problem which naturally provides us with the necessary and sufficient information, and we comment, among other issues, on a dispensable third part of the assumption. 
\end{abstract}

\numberwithin{equation}{section}
\section{Introduction}\label{intro}

The study of the fundamental mathematical model for the kinetics of the general enzymatic reaction with chemical equation
\begin{equation}
\label{SEP}
S+E\rightarrow E+P,
\end{equation}
where $S$ is the substrate, $P$ is the product and  $E$ is the enzyme that catalyses it, has a long history, which we briefly present below.

Already since 1894, Fischer \cite{fischer1894einfluss} derived the \textit{lock and key model} for the interpretation of biocatalysis. Already prior to 1901, Brown suggested an intermediate step in the enzymatic reaction that is described by (\ref{SEP}), where the substrate forms a \textit{complex} with the enzyme before the beginning of the catalysis, an idea that was eventually published in 1902 \cite{brown1902}. Thus, it had already been realised by that time, that enzymatic biochemical reactions should take place in at least two stages, and in fact these stages should have different  {\textit time scales}. Based on this idea, combined with conversations he had with Bodenstein, Henri published in 1902 \cite{henri1902gesetz} and then in 1903 \cite{henri1903} an initial version of a reliable differential equation for the description of the kinetics of the enzyme reaction with chemical equation given by (\ref{SEP}), an idea he had conceived as early as 1901. A decade later, in 1913, Michaelis and Menten \cite{michment} (translated in English in \cite{michmenten}), extracted this equation by using a more detailed and analytical form that makes use of the \textit{rapid equilibrium assumption}; they interpreted it convincingly and studied it thoroughly. In particular, using as an example the invertase-catalysed hydrolysis of sucrose into glucose and fructose, they studied (\ref{SEP}) through the chemical mechanism 
\begin{equation}
\label{SECP}
S+E\,\underset{k_{-1}}{\overset{k_1}{\rightleftharpoons}}\,C\xrightarrow{k_2} E+P,
\end{equation}
where $k_1,k_{-1},k_2\,>0$ and $C$ represents the substrate-enzyme complex, and indirectly concluded that, when
\begin{equation}
\label{appSE}
\left[S\right]={\left[S\right]}_0\gg{\left[E\right]}_0=\left[E\right]\text{ and }\left[C\right]={\left[C\right]}_0=0,\text{ for }t=0,
\end{equation}
a condition acceptable in enzymatic reactions, then for the rate $\upsilon$ of the enzymatic reaction with chemical equation (\ref{SEP}) it holds that
\begin{equation}
\label{vappMM}
\upsilon\approx\frac{\upsilon_{\sup}\left[S\right]}{K_{{\rm dis}}+\left[S\right]}\,,
\end{equation}
where $$\upsilon_{\sup}\coloneqq k_2{\left[E\right]}_0\,,$$ and
\begin{equation}
\label{KSconst}
K_{{\rm dis}}\coloneqq \frac{k_{-1}}{k_1\,}. 
\end{equation}
$K_{{\rm dis}}$ is the constant that is nowadays called the \textit{dissociation constant} (of the complex).

On the other hand, Van Slyke and Cullen, working in parallel with  Michaelis and Menten, but studying urease-catalysed hydrolysis of urea to ammonia and carbon dioxide, used - instead of (\ref{SECP}) - the chemical mechanism $$S+E \xrightarrow{k_1} C\xrightarrow{k_2} E+P$$ and concluded in 1914 \cite{vanslykecullen} to
\begin{equation}
\label{vappVSC}
\upsilon\approx\frac{\upsilon_{\sup}\left[S\right]}{K_{\!V\!SC}+\left[S\right]},
\end{equation}
instead of (\ref{vappMM}), where 
\begin{equation}
\label{KVSCconst}
K_{\!V\!SC}\coloneqq\frac{k_{2}}{k_1}\,.
\end{equation}
$K_{\!V\!SC}$ is a constant that is now known as the \textit{Van Slyke-Cullen constant}.  

In 1925, Briggs and Haldane \cite{briggs1925note} published a short note where they composed the ideas of Michaelis \& Menten and  Van Slyke \& Cullen through a raw first version of a new assumption, known presently as the \textit{standard quasi-steady-state assumption}. In particular, they improved (\ref{vappMM}) and (\ref{vappVSC}), demonstrating that
\begin{equation}
\label{MMeq}
\upsilon\approx\frac{\upsilon_{\sup}\left[S\right]}{K_{\!M}+\left[S\right]}\,,
\end{equation}
where 
\begin{equation}
\label{KMconst}
K_{\!M}\coloneqq K_{{\rm dis}}+K_{\!V\!SC} =\frac{k_{-1}+k_2}{k_1}\,.
\end{equation}
It is a standard expression nowadays, that \textquote{(\ref{MMeq}) characterises the Michaelis-Menten kinetics}, and the constant (\ref{KMconst}) is called the  Michaelis-Menten constant.

Lineweaver and Burk in 1934 \cite{lineweaver1934determination} established (\ref{MMeq}) in the form $$\frac{1}{\upsilon}\approx \frac{1}{\upsilon_{\sup}}+\frac{K_{\!M}}{\upsilon_{\sup}}\frac{1}{\left[S\right]}\,,$$  as a tool for experimental calculation of the values $\upsilon_{\sup}$ and $K_{\!M}$. 

Already since the beginning of the second half of the 20th century and throughout it, many researchers have dealt with the validity of the quasi-steady-state assumption and the determination of the two time scales of the model through the application of  {\it perturbation methods}. However, it was much later, in 1988 and in 1989, when Segel \cite{segel1988validity}, and Segel and  Slemrod \cite{segel1989quasi}, respectively, showed that if
\begin{equation}
\label{appSEK}
\left(K_{\!M} +{\left[S\right]}_0\right)\gg{\left[E\right]}_0\text{ and }{\left[C\right]}_0=0,
\end{equation}
then there are indeed two time scales, which recorded as follows
\begin{equation*}
t_C=\frac{1}{k_1\left(K_{\!M} +{\left[S\right]}_0\right)}\ll\frac{K_{\!M} +{\left[S\right]}_0}{k_2{\left[E\right]}_0}=t_S,
\end{equation*}
and it holds that
\begin{equation}
\label{upslnsss}
\upsilon\approx\begin{cases}
0,&\text{ for times comparable to }t_C,\\
\text{as in (\ref{MMeq})},&\text{ for times comparable to }t_S.
\end{cases}
\end{equation}
In fact, (\ref{appSEK}) is more general than (\ref{appSE}) since it allows  $$\frac{{\left[E\right]}_0}{{\left[S\right]}_0}=O{\left(1\right)},\text{ as }\frac{{\left[E\right]}_0}{K_{\!M}+{\left[S\right]}_0}\rightarrow 0^+,$$ or even $${\left[E\right]}_0\gg {\left[S\right]}_0.$$ 

In 1997, Schnell and Mendoza \cite{schnell1997closed} captured the solution of the Michaelis-Menten kinetics equation in closed form, by using the Lambert $W$ function, and in particular a special case of it, which is defined by its inverse  as follows $$W^{-1}{\left(x\right)}=x\,\exp{\left\{x\right\}},\text{ }\forall x\geq 0.$$ In addition, via the aforementioned work of Segel and Slemrod, an initial form of another hypothesis was introduced for the first time,  the \textit{reverse quasi-steady-state assumption}, and it was shown that, when
\begin{equation}
\label{SSrA}
{\left[E\right]}_0\gg K_{\!M}\text{ and }{\left[C\right]}_0=0,
\end{equation}
then there are again two time scales, $$\tilde{t}_S=\frac{1}{k_1{\left[E\right]}_0}\ll\frac{1}{k_2}=\tilde{t}_C,$$ and it holds that 
\begin{equation}
\label{upslnssr}
\upsilon\approx\begin{cases}
k_1{\left[E\right]}_0\left[S\right],&\text{ for times comparable to }\tilde{t}_S,\\
0,&\text{ for times comparable to }\tilde{t}_C.
\end{cases}
\end{equation}

About a decade later, in 2000, Schnell and Maini \cite{schnell2000enzyme} found that (\ref{SSrA}) is not sufficient for (\ref{upslnssr}) to hold; on the contrary, the new case should have the form
\begin{equation}
\label{SSrA2}
{\left[E\right]}_0\gg K_{\!M},\text{ }{\left[E\right]}_0\gg K_{\!M}\text{ and }{\left[C\right]}_0=0.
\end{equation}

Finally, let us mention that in 1996, with the work of  Borghans, Boer and Segel \cite{borghans1996extending}, the total substrate concentration, $\left[T\right]$ is introduced, i.e., the sum of the concentration of the unbound/free substrate plus the concentration of the bound substrate in the form of complex with the enzyme, that is
\begin{equation}
\label{TclneqSC}
\left[T\right]\coloneqq \left[S\right]+\left[C\right],
\end{equation}
to describe an alleged third hypothesis that shares common ground with both the previous ones, the so called \textit{total quasi-steady-state assumption}, and since then several researchers have adopted, and dealt with this hypothesis.\\

In this work, our novel results are
\begin{enumerate}
\item We propose a general and simple algorithm for the proper scaling of every problem with non negative solutions in a bounded domain, and we essentially employ it in our case. Until now, only a ``rough'' rule for the non dimensionalisation process is utilised in applications, which states that the scales considered for the variables of a problem are chosen so that they should be roughly of the same order of magnitude of the respective variables themselves \cite{logan2013applied}. The proposed procedure is as follows: 
\begin{enumerate}[label=\Roman*.]
\item Identification of the bounded feasible region - i.e., the range of the dependent variables - of the problem. 
\item Scaling of the dependent variables of the problem by their respective supremum feasible values - which do exist, since the non negative solutions exist in a bounded domain. 
\item Natural scaling of the independent variables by gathering the remaining terms of the previous step. 
\end{enumerate}
By this algorithm, 
\begin{itemize}
\item[$\circ$] the dependent variables are comparable with each other, since they all range onto $\left[0,1\right]$,  
\item[$\circ$] any scale of the independent variables follows naturally by the process, hence there is no need of the - somehow unjustified - approach of \textquote{an estimate of the minimum value for which the variable undergoes a significant change in magnitude} (see, e.g., \cite{lin1988mathematics}, \cite{segel1988validity} and \cite{segel1989quasi}), which is widely adopted thenceforth, for the choice of the largest of the two time scales (in our case),
\item[$\circ$] the quantity $\varepsilon$, that characterises both the standard and the reverse quasi-steady-state assumptions, arises effortless by the problem itself.
\end{itemize}
\item We clarify the purely quantitative nature of the standard and the reverse quasi-steady-state assumptions, as opposed to the aforementioned qualitative one. In particular, (\ref{sQSSA}) and (\ref{rQSSA}) do not serve for the validation of the standard and the reverse, respectively, quasi-steady-state assumptions: they define them. 
\item We relinquish the, so called, {\it total quasi-steady-state assumption}, by showing that there is no substantive third hypothesis, but only a different approach to the first two. We note that such a duality, characterised by a positive parameter $\varepsilon$, that either tends to $0$ or to $\infty$, is common in applications, for instance in the study of Hamiltonian systems possessing either a relatively small or a relatively large Hamiltonian. 
\end{enumerate}

As far as the mathematical tools employed in the present work are concerned, we note the following: 
\begin{itemize}
\item For the sake of brevity, we neither state nor discuss the necessary concepts and fundamental results regarding a Cauchy (initial value) problem for vector first order autonomous, or not, nonlinear ordinary differential equations. These issues would certainly comprise the existence, the uniqueness, the extendibility, the regularity, the continuous and smooth dependence of the solutions on the data (initial condition, vector field determining the differential equation, possible parameters appearing in the Cauchy problem). Additional principal issues would include notions of stability, local and global techniques for studying it, and essential results of the qualitative theory of ODEs, in general. There is a huge literature on these topics; indicatively, we refer to   \cite{haleodes}, \cite{meissdiffdynsyst}, \cite{pontryaginodes},  and \cite{strogatz2018nonlinear}.
\item A powerful technique for problems with relatively small (or large) parameters is the {\it Method of Matched Asymptotics}, where approximate solutions, accurate in one region of the problem domain, are matched to different approximate solutions, accurate in another region. This subject is discussed in many books, see, e.g., \cite{bender2013advanced}, \cite{holmes2012pertmeth}, \cite{lin1988mathematics},  and \cite{logan2013applied}. 
\end{itemize}

\section{Principal analysis of the problem} 

In this section, we introduce the main problem and proceed to its basic analysis, that comprises the identification of the feasible regions, the well posedness of the problem, the determination of the invariant sets, as well as  the simplification and the stability analysis of the problem. 

\subsection{Cauchy problem}

Employing the chemical mechanism (\ref{SECP}) along with the Law of Mass Action \cite{voit2015150}, we arrive at the equations
\begin{subequations}
\label{SECPeq}
\begin{align}
\dfrac{\de \left[S\right]}{\de t}&=-k_1\left[S\right]\left[E\right]+k_{-1}\left[C\right]\label{SECPeq;a},\\
\dfrac{\de \left[E\right]}{\de t}&=-k_1\left[S\right]\left[E\right]+\left(k_{-1}+k_2\right)\left[C\right]\label{SECPeq;b},\\
\dfrac{\de \left[C\right]}{\de t}&=k_1\left[S\right]\left[E\right]-\left(k_{-1}+k_2\right)\left[C\right]\label{SECPeq;c},\\
\dfrac{\de \left[P\right]}{\de t}&=k_2\left[C\right],\label{SECPeq;d}
\end{align}
\end{subequations}
and the corresponding Cauchy problem reads:
\[\label{SECP_Cauchy}\tag{$SECP$}
\begin{minipage}{.8\textwidth}
\centering 
\textit{Given ${\left[S\right]}_0,{\left[E\right]}_0,{\left[C\right]}_0,{\left[P\right]}_0\,\geq 0$, we seek an interval $\mathcal{I}\subseteq\mathbb{R}$ with $0\in\mathcal{I}$, and a function $\left(\left[S\right],\left[E\right],\left[C\right],\left[P\right]\right)\colon\mathcal{I}\rightarrow{\left[0,\infty\right)}^4$, such that $\left(\left[S\right],\left[E\right],\left[C\right],\left[P\right]\right)$ satisfies both \eqref{SECPeq} in ${\left(\mathcal{I}\setminus\left\{0\right\}\right)}^\circ$ and $\left(\left[S\right],\left[E\right],\left[C\right],\left[P\right]\right)=\left({\left[S\right]}_0,{\left[E\right]}_0,{\left[C\right]}_0,{\left[P\right]}_0\right)$, for $t=0$.}
\end{minipage}
\]

For a solution of (\ref{SECP_Cauchy}) it holds that $$\dfrac{\de \left[S\right]}{\de t}+\dfrac{\de \left[C\right]}{\de t}+\dfrac{\de \left[P\right]}{\de t}=0,$$ or, equivalently, 
\begin{equation}
\label{SplsCplsP}
\left[S\right]+\left[C\right]+\left[P\right]={\left[S\right]}_0+{\left[C\right]}_0+{\left[P\right]}_0\eqqcolon A_1,
\end{equation}
due to the initial condition of (\ref{SECPeq}), as well as that $$\dfrac{\de \left[E\right]}{\de t}+\dfrac{\de \left[C\right]}{\de t}=0,$$ or, equivalently, 
\begin{equation}
\label{EplsC}
\left[E\right]+\left[C\right]={\left[E\right]}_0+{\left[C\right]}_0\eqqcolon A_2.
\end{equation}
From (\ref{SplsCplsP}) and (\ref{EplsC}) combined with the non negativity of the components of the solutions of (\ref{SECP_Cauchy}), we conclude that
\begin{equation}
\label{SCEPA1A2v1}
\left[S\right]\leq A_1,\text{ }\left[E\right]\leq A_2,\text{ }\left[C\right]\leq \min{\left\{A_1,A_2\right\}}\text{ and }\left[P\right]\leq A_1.
\end{equation}
In addition, from  (\ref{SECPeq;c}) together with the bounds for $\left[S\right]$ and $\left[E\right]$ in (\ref{SCEPA1A2v1}) we have that
\begin{equation}
\label{Canthrinq}
\left[C\right]\leq\frac{A_1A_2}{K_{\!M}},
\end{equation}
where $K_{\!M}$ is defined as in (\ref{KMconst}), whereas the rest of the equations of (\ref{SECPeq}) do not include further related information. Thus, from  (\ref{SCEPA1A2v1}) and (\ref{Canthrinq}) we finally get that
\begin{equation}
\label{SCEPA1A2}
\left[S\right]\leq A_1,\text{ }\left[E\right]\leq A_2,\text{ }\left[C\right]\leq \min{\left\{A_1,A_2,\frac{A_1A_2}{K_{\!M}}\right\}}\eqqcolon A_3\text{ and }\left[P\right]\leq A_1.
\end{equation}

In the light of (\ref{SCEPA1A2}), we set
$$\Omega_0\coloneqq\left\{\left(s,e,c,p\right)\in\left[0,A_1\right]\times\left[0,A_2\right]\times\left[0,A_3\right]\times\left[0,A_1\right]\,\big|\,s+c+p=A_1,\text{ }e+c=A_2\right\}$$ and we can therefore consider an equivalent to (\ref{SECP_Cauchy})  problem as follows:
\[
\begin{minipage}{.8\textwidth}
\centering 
\textit{Given ${\left[S\right]}_0,{\left[E\right]}_0,{\left[C\right]}_0,{\left[P\right]}_0\,\geq 0$,  we are looking for an interval $\mathcal{I}\subseteq\mathbb{R}$ with $0\in\mathcal{I}$ and a function $\left(\left[S\right],\left[E\right],\left[C\right],\left[P\right]\right)\colon\mathcal{I}\rightarrow \Omega_0$, such that $\left(\left[S\right],\left[E\right],\left[C\right],\left[P\right]\right)$ to satisfy both \eqref{SECPeq} in ${\left(\mathcal{I}\setminus\left\{0\right\}\right)}^\circ$ and $\left(\left[S\right],\left[E\right],\left[C\right],\left[P\right]\right)=\left({\left[S\right]}_0,{\left[E\right]}_0,{\left[C\right]}_0,{\left[P\right]}_0\right)$ for $t=0$.}
\end{minipage}
\]

Employing standard arguments of the theory of Ordinary Differential Equations,  
we can conclude that (\ref{SECP_Cauchy}) is globally well posed, with an infinitely smooth solution in an interval $\widetilde{\mathcal{I}}\subseteq\mathbb{R}$, where $$\widetilde{\mathcal{I}}=\mathbb{R},\text{ or }\widetilde{\mathcal{I}}=\left[-a,\infty\right)\,\text{for some }a\in\left[0,\infty\right).$$ In addition, when $A_2=0$, the unique solution  is the constant  $$\left(S,E,C,P\right)=\left({\left[S\right]}_0,0,0,{\left[P\right]}_0\right),\text{ }\forall t\in\mathbb{R}.$$ Thus, when $A_2=0$, $\Omega_0$ reduces to $$\Omega_0=\left\{\left(s,0,0,p\right)\in\left[0,A_1\right]\times{\left\{0\right\}}^2\times\left[0,A_1\right]\,\big|\,s+p=A_1\right\}$$ which is invariant (in particular, every singleton $\left\{\left(s,0,0,A_1-s\right)\right\}$ for $s\in\left[0,A_1\right]$ is invariant), whereas $\Omega_0$ is positively invariant, when $A_2>0$. 

\subsection{A simpler equivalent problem}

Given that (\ref{SplsCplsP}) and (\ref{EplsC}) hold, we conclude that system (\ref{SECPeq})  can be equivalently reduced to
\begin{subequations}
\label{SCeq}
\begin{align}
\dfrac{\de \left[S\right]}{\de t}&=-k_1A_2\left[S\right]+k_1\left[S\right]\left[C\right]+k_{-1}\left[C\right],\label{SCeq;a}\\
\dfrac{\de \left[C\right]}{\de t}&=k_1A_2\left[S\right]-k_1\left[S\right]\left[C\right]-\left(k_{-1}+k_2\right)\left[C\right]\label{SCeq;b}.
\end{align}
\end{subequations}

Let us now study the above subsystem. Using (\ref{SCeq;b})  combined with the bound of $\left[S\right]$ in (\ref{SCEPA1A2}) we have that
\begin{equation}
\label{CA12}
\left[C\right]\leq \frac{A_1A_2}{K_{\!M}+A_1}\,.
\end{equation}
Therefore, from the bound of  $\left[C\right]$ in (\ref{SCEPA1A2}) and from (\ref{CA12}) we eventually get that
\begin{equation}
\label{minC}
\left[C\right]\leq \min{\left\{A_1,A_2,\frac{A_1A_2}{K_{\!M}},\frac{A_1A_2}{K_{\!M}+A_1}\right\}}=\min{\left\{A_1,\frac{A_1A_2}{K_{\!M}+A_1}\right\}}\eqqcolon A_4.
\end{equation}
In fact, key to what follows are the immediately verifiable inferences
\begin{equation}
\label{bscA12}
\fbox{$A_2\leq K_{\!M}+A_1\Rightarrow A_4=\dfrac{A_1A_2}{K_{\!M}+A_1}$},
\end{equation}
and on the other hand, 
\begin{equation}
\label{bscA12neg}
\fbox{$A_2\geq K_{\!M}+A_1\Rightarrow A_4=A_1$}. 
\end{equation}
Now, in the light of the bound for $\left[S\right]$ in (\ref{SCEPA1A2}) and of (\ref{minC}), we set $$\Omega_1\coloneqq\left\{\left(s,c\right)\in\left[0,A_1\right]\times\left[0,A_4\right]\,\big|\,s+c\leq A_1\right\}$$ and so we can consider the equivalent, to (\ref{SECP_Cauchy}),  problem as follows:
\[\label{SC_Cauchy}\tag{$SC$}
\begin{minipage}{.8\textwidth}
\centering 
\textit{Given ${\left[S\right]}_0,{\left[E\right]}_0,{\left[C\right]}_0,{\left[P\right]}_0\,\geq 0$, we seek an interval $\mathcal{I}\subseteq\mathbb{R}$ with $0\in\mathcal{I}$, and a function $\left(\left[S\right],\left[C\right]\right)\colon\mathcal{I}\rightarrow \Omega_1$, such that $\left(\left[S\right],\left[C\right]\right)$ satisfies both \eqref{SCeq} in ${\left(\mathcal{I}\setminus\left\{0\right\}\right)}^\circ$, and $\left(\left[S\right],\left[C\right]\right)=\left({\left[S\right]}_0,{\left[C\right]}_0\right)$ for $t=0$.}
\end{minipage}
\]

\subsection{Stability analysis}

First, we can easily deduce that $$\begin{cases}
\left(s,0,0,A_1-s\right)\text{ where }s\in\left[0,A_1\right],&\text{ when }A_2=0,\\
\left(0,e,0,A_1\right)\text{ where }e\in\left[0,A_2\right],&\text{ when }A_2>0,
\end{cases}$$
are the steady states of  (\ref{SECP_Cauchy}). However, we immediately conclude that it makes sense to study their stability only for the non-trivial case, where  $$A_1>0\text{ and }A_2>0.$$ It is sufficient though, as usually, to study the stability of $\left(0,0\right)$, as a steady state of (\ref{SC_Cauchy}), when $A_1>0$ and $A_2>0$.

As for the local stability of $\left(0,0\right)$, we calculate the Jacobi matrix: 
\begin{equation*}
J{\left(s,c\right)}=
\begin{pmatrix}
k_1\left(c-A_2\right) & k_1 s+k_{-1} \\
k_1\left(A_2-c\right) & -k_1 s-\left(k_{-1}+k_2\right) 
\end{pmatrix}.
\end{equation*}
Its eigenvalues  at  $\left(0,0\right)$, are
$$\lambda_{\pm}=\frac{1}{2}\left(-k_1 A_2-\left(k_{-1}+k_2\right)\pm{\left({\left(k_1 A_2+\left(k_{-1}+k_2\right)\right)}^2-4k_1k_2A_2\right)}^{\frac{1}{2}}\right).$$ Since $${\left(k_1 A_2+\left(k_{-1}+k_2\right)\right)}^2-4k_1k_2A_2={\left(k_1 A_2+\left(k_{-1}-k_2\right)\right)}^2+4k_{-1}k_2\geq 0,$$  the origin $\left(0,0\right)$ is locally asymptotically stable for (\ref{SC_Cauchy}), since $$\lambda_{\pm}<0.$$ In fact, we can also find, as usually, a local approach to the solution close to  $\left(0,0\right)$; we omit it for the sake of brevity.

As for the global stability of $\left(0,0\right)$, we can apply the Bendixson-Dulac Negative Criterion with 
\begin{align*}
g\,\colon{\left(\mathbb{R}_+\right)}^2&\rightarrow \left(0,\infty\right)\\
\left(s,c\right)&\mapsto g{\left(s,c\right)}=\dfrac{1}{sc},
\end{align*}
thereby obtaining the desired result, since in ${\left(\Omega_1\right)}^\circ$ it holds that
\begin{equation*}
\mathrm{div}{\left[\left(\frac{-k_1A_2}{\left[C\right]}+k_1+\frac{k_{-1}}{\left[S\right]},\frac{k_1A_2}{\left[C\right]}-k_1-\frac{\left(k_{-1}+k_2\right)}{\left[S\right]}\right)\right]}=-\frac{k_{-1}}{{\left[S\right]}^2}-\frac{k_1A_2}{{\left[C\right]}^2}<0. 
\end{equation*}

\section{The standard quasi-steady-state assumption}

The standard quasi-steady-state assumption is 
\begin{equation}
\label{sQSSA}\tag{$sQSSA$}
\fbox{$A_1>0$ and $0<A_2\ll K_{\!M}+A_1$},
\end{equation}
or, equivalently, $$A_1>0\text{ and }0<A_2\ll K_{\!M},\text{ or }0<A_2\ll A_1\text{ (or these two together)},$$ and provided that it holds, we study (\ref{SC_Cauchy}). 

We consider two approaches for examining the assumption, the free substrate approach, where the concentration dynamics of the unbound substrate, $\left[S\right]$, is studied, and the total substrate approach,  where the concentration dynamics of the total substrate, $\left[T\right]$, is studied, as defined in (\ref{TclneqSC}). 

\subsection{Free substrate approach}

Using only (\ref{sQSSA}) we will show that
\begin{enumerate}
\item Problem (\ref{SC_Cauchy}), and therefore problem (\ref{SECP_Cauchy}) as well, has inherently two time scales which we will determine. In fact, (\ref{sQSSA}) owes its name to the existence of the above time scales. In particular, except for a short initial time interval, where the enzymatic reaction with chemical equation  (\ref{SEP})  is not evolving, i.e., $\upsilon\approx 0$, during the rest of the time the enzymatic reaction is at a \textquote{steady state}, in which (\ref{MMeq}) holds.
\item There is a good uniform approximation in closed form to the solution of (\ref{SC_Cauchy}), and therefore to  (\ref{SC_Cauchy}) as well,  which we will determine.
\end{enumerate}
To highlight the above time scales, the first and basic step is scaling (\ref{SC_Cauchy}). Thus, as usually, due to the bound of $\left[S\right]$ in (\ref{SCEPA1A2}) and the relation (\ref{minC}), we choose the dimensionless dependent variables as 
 $$S_\alpha{\left(t_\alpha\right)}\coloneqq\frac{1}{A_1}\left[S\right]{\left(\frac{t}{t_*}\right)}\text{ and }C_\alpha{\left(t_\alpha\right)}\coloneqq\frac{1}{A_4}\left[C\right]{\left(\frac{t}{t_*}\right)},$$ where we have chosen an arbitrary, for the time being, time scale $t_*>0$  for the scaling, i.e., $$t_\alpha\coloneqq\frac{t}{t_*},$$ the determination  of which will arise in a natural manner during the process. We note, however, that - given (\ref{sQSSA}) - it follows from (\ref{bscA12}) that 
$$\fbox{$A_4=\varepsilon A_1\ll A_1$}$$ and $$C_\alpha{\left(t_\alpha\right)}=\frac{1}{\varepsilon A_1}\left[C\right]{\left(\frac{t}{t_*}\right)},$$ where
\begin{equation}
\label{epsln}
0\,\overset{\text{(\ref{sQSSA})}}{<}\,\varepsilon\coloneqq\frac{A_2}{K_{\!M}+A_1}\,\overset{\text{(\ref{sQSSA})}}{\ll}\,1, 
\end{equation}
i.e., equivalently, 
\begin{equation}
\label{epsln2}
0<\varepsilon\approx 0. 
\end{equation}
A first conclusion is that the possible change of  $\left[S\right]$ is much larger than the corresponding one of  $\left[C\right]$. Now, (\ref{SCeq}) will take the following form
\begin{subequations}
\label{SCeqdmlss}
\begin{align}
\dfrac{\de S_\alpha}{\de t_\alpha}&=k_1A_2t_*\left(-S_\alpha+\frac{\sigma}{1+\sigma}S_\alpha C_\alpha+\frac{\rho}{\left(1+\rho\right)\left(1+\sigma\right)}C_\alpha\right)\label{SCeqdmlss;a},\\
\dfrac{\de C_\alpha}{\de t_\alpha}&=k_1\left(K_{\!M}+A_1\right)t_*\left(S_\alpha-\frac{\sigma}{1+\sigma}S_\alpha C_\alpha-\frac{1}{1+\sigma}C_\alpha\right)\label{SCeqdmlss;b}, 
\end{align}
\end{subequations}
where 
\begin{equation}
\label{sgmrh}
\sigma\coloneqq\frac{A_1}{K_{\!M}}\text{ and }\rho\coloneqq\frac{K_{\rm dis}}{K_{\!V\!SC}}=\frac{k_{-1}}{k_2},
\end{equation}
where $K_{\rm dis}$ and $K_{\!V\!SC}$ are as in (\ref{KSconst}) and (\ref{KVSCconst}), respectively. 

Observing (\ref{SCeqdmlss}), we define 
\begin{equation}
\label{tCtS}
\frac{1}{k_1\left(K_{\!M}+A_1\right)}\eqqcolon t_1\,\overset{\text{(\ref{epsln})}}{\ll}\, \frac{t_1}{\varepsilon}=\frac{1}{k_1 A_2}\eqqcolon t_2, 
\end{equation}
to conclude that
\begin{equation*}
t_*=t_1\Rightarrow 
\begin{cases}
k_1 A_2 t_*=\varepsilon,\\
k_1\left(K_{\!M}+A_1\right)t_*=1,
\end{cases}
\text{ and }\,\,\,t_*=t_2\Rightarrow 
\begin{cases}
k_1 A_2 t_*=1,\\
k_1\left(K_{\!M}+A_1\right)t_*=\frac{1}{\varepsilon},
\end{cases}
\end{equation*}
and so (\ref{SCeqdmlss}) takes the following form:
\begin{itemize}
\item If \fbox{$t_\alpha=\dfrac{t}{t_1}$}, then 
\begin{subequations}
\label{SCeqdmlsst1}
\begin{align}
\dfrac{\de S_\alpha}{\de t_\alpha}&=\varepsilon\left(-S_\alpha+\frac{\sigma}{1+\sigma}S_\alpha C_\alpha+\frac{\rho}{\left(1+\rho\right)\left(1+\sigma\right)}C_\alpha\right)\label{SCeqdmlsst1;a},\\
\dfrac{\de C_\alpha}{\de t_\alpha}&=S_\alpha-\frac{\sigma}{1+\sigma}S_\alpha C_\alpha-\frac{1}{1+\sigma}C_\alpha\label{SCeqdmlsst1;b}. 
\end{align}
\end{subequations}
\item If \fbox{$t_\alpha=\dfrac{t}{t_2}$}, then 
\begin{subequations}
\label{SCeqdmlsst2}
\begin{align}
\dfrac{\de S_\alpha}{\de t_\alpha}&=-S_\alpha+\frac{\sigma}{1+\sigma}S_\alpha C_\alpha+\frac{\rho}{\left(1+\rho\right)\left(1+\sigma\right)}C_\alpha\label{SCeqdmlsst2;a},\\
\dfrac{\de C_\alpha}{\de t_\alpha}&=\frac{1}{\varepsilon}\left(S_\alpha-\frac{\sigma}{1+\sigma}S_\alpha C_\alpha-\frac{1}{1+\sigma}C_\alpha\right)\label{SCeqdmlsst2;b}. 
\end{align}
\end{subequations}
\end{itemize}

Setting $$\Omega_\varepsilon\coloneqq\left\{\left(s,c\right)\in{\left[0,1\right]}^2\,\big|\,s+\varepsilon c\leq 1\right\},$$ the scaled version of (\ref{SC_Cauchy}) is: 
\[\label{SCdmlss_Cauchy}\tag{$SC\alpha_s$}
\begin{minipage}{.8\textwidth}
\centering 
\textit{Given ${\left[S\right]}_0,{\left[C\right]}_0,{\left[P\right]}_0\,\geq 0$ and $\varepsilon>0$, we seek an interval $\mathcal{I}\subseteq\mathbb{R}$ with $0\in\mathcal{I}$, and a function $\left(S_\alpha,C_\alpha\right)\colon\mathcal{I}\rightarrow \Omega_\varepsilon$, such that $\left(S_\alpha,C_\alpha\right)$ satisfying both \eqref{SCeqdmlsst1} if $t_\alpha=\dfrac{t}{t_1}$ or \eqref{SCeqdmlsst2} if $t_\alpha=\dfrac{t}{t_2}$ in ${\left(\mathcal{I}\setminus\left\{0\right\}\right)}^\circ$ and $\left(S_\alpha,C_\alpha\right)=\left(\dfrac{{\left[S\right]}_0}{A_1},\dfrac{{\left[C\right]}_0}{\varepsilon A_1}\right)$, for $t_\alpha=0$.}
\end{minipage}
\]
\begin{figure}[!htbp]
\centering
\includegraphics[width=.5\linewidth]{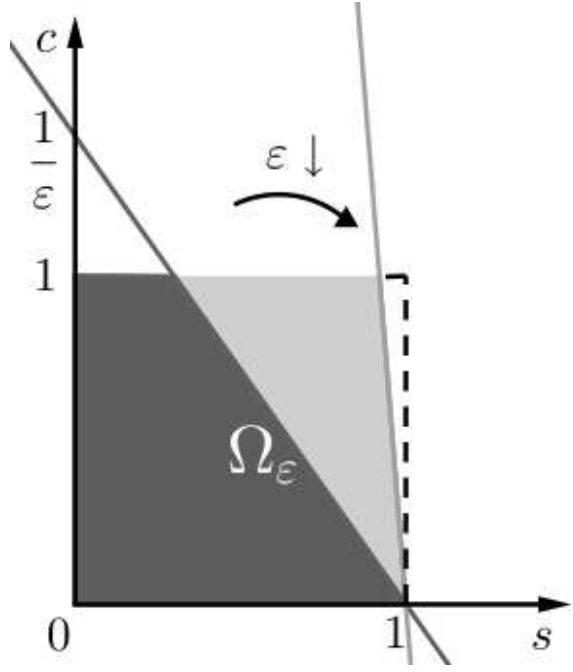}
\caption{The invariant set $\Omega_\varepsilon$ of problem \eqref{SCdmlss_Cauchy}. We notice that $\Omega_\varepsilon\rightarrow{\left[0,1\right]}^2$ as $\varepsilon\rightarrow 0^+$.}
\label{Omega-epsilon}
\end{figure}

We study separately each of the two versions of (\ref{SCdmlss_Cauchy}) in order to find an inner and outer, respectively, approximation to the solution $\left(\left[S\right],\left[C\right]\right)$ of (\ref{SC_Cauchy}), i.e., one approximation for times  \textit{comparable} to $t_1$ and another one for  times  \textit{comparable} to $t_2$, respectively. In more detail: 
\begin{itemize}
\item Looking at (\ref{SC_Cauchy}) as a \textit{perturbed} problem, with perturbation  $\varepsilon>0$ close to $0$, we have the following information on (\ref{SCeqdmlsst1;a}) $$\dfrac{\de S_\alpha}{\de t_\alpha}\,=\,\underbrace{\varepsilon\left(\underbrace{-S_\alpha+\frac{\sigma}{1+\sigma}S_\alpha C_\alpha+\frac{\rho}{\left(1+\rho\right)\left(1+\sigma\right)}C_\alpha}_{\substack{=\,O{\left(1\right)}\text{ uniformly (with respect to $t_\alpha$) as }\varepsilon\rightarrow 0^+,\\ \text{since }\left(S_\alpha,C_\alpha,\right)\,=\,O{\left(1\right)}\text{ uniformly as }\varepsilon\rightarrow 0^+,\\ \text{since }\left(S_\alpha,C_\alpha\right)\,\in\, \Omega_\varepsilon\text{ }\forall\varepsilon>0}}\right)}_{=\,O{\left(\varepsilon\right)}\text{ uniformly as }\varepsilon\rightarrow 0^+},$$ as well as $$\underbrace{\dfrac{\de S_\alpha}{\de t_\alpha}}_{\substack{=\,O{\left(1\right)}\text{ uniformly as }\varepsilon\rightarrow 0^+,\\ \text{when }t_\alpha\,=\,O{\left(1\right)}\text{ uniformly as }\varepsilon\rightarrow 0^+,\\ \text{since }S_\alpha\,=\,O{\left(1\right)}\text{ uniformly as }\varepsilon\rightarrow 0^+\\ \text{combined with the definition of the derivative}}}\,=\,\underbrace{\dots}_{=\,O{\left(\varepsilon\right)}\text{ uniformly as }\varepsilon\rightarrow 0^+},$$ thus 
$$\underbrace{\frac{\de S_\alpha}{\de t_\alpha}}_{\substack{=\,O{\left(\varepsilon\right)}\text{ uniformly as }\varepsilon\rightarrow 0^+,\\ \text{when }t_\alpha\,=\,O{\left(1\right)}\text{ uniformly as }\varepsilon\rightarrow 0^+}}.$$ Thus, due to (\ref{epsln2}) it follows that
\begin{equation}
\label{deSalphat1approx}
\frac{\de S_\alpha}{\de t_\alpha}\approx 0,\text{ when }\exists B>0\text{ independent of }\varepsilon,\text{ such that: }\left|t_\alpha\right|\leq B,
\end{equation}
and due to the initial condition of (\ref{SCdmlss_Cauchy}) we eventually have that $$S_\alpha\approx \frac{{\left[S\right]}_0}{A_1},\text{ when }\exists B>0:\text{ }\left|t_\alpha\right|\leq B.$$ If we insert the above approximate equality in (\ref{SCeqdmlsst1;b}), then the later becomes an approximate linear differential equation, the solution of which is
\begin{equation*}
C_\alpha\approx \frac{\left(1+\sigma\right)S_\alpha}{1+\sigma S_\alpha}+\left(\dfrac{{\left[C\right]}_0}{\varepsilon A_1}-\frac{\left(1+\sigma\right)S_\alpha}{1+\sigma S_\alpha}\right)\exp{\left\{-\frac{1+\sigma S_\alpha}{1+\sigma}\, t_\alpha\right\}},\text{ when }\exists B>0:\text{ }\left|t_\alpha\right|\leq B,
\end{equation*}
given the initial condition of (\ref{SCdmlss_Cauchy}).\\ Therefore, the inner approximation, $\left({\left[S\right]}_{{\rm in}},{\left[C\right]}_{{\rm in}}\right)$, of the solution $\left(\left[S\right],\left[C\right]\right)$ of (\ref{SC_Cauchy}), i.e., the approximation for those $t$ for which it holds that $$\text{\fbox{$\exists B>0$ independent of $\varepsilon$, such that: $\left|t\right|\leq B t_1$}},$$ is $$\left({\left[S\right]}_{{\rm in}},{\left[C\right]}_{{\rm in}}\right)=\left({\left[S\right]}_0,\frac{A_2{\left[S\right]}_0}{K_{\!M}+{\left[S\right]}_0}+\left({\left[C\right]}_0-\frac{A_2{\left[S\right]}_0}{K_{\!M}+{\left[S\right]}_0}\right)\exp{\left\{-k_1\left(K_{\!M}+{\left[S\right]}_0\right)t  \right\}}\right).$$
\item For (\ref{SCeqdmlsst2;b}) we have that $$\dfrac{\de C_\alpha}{\de t_\alpha}=\underbrace{\frac{1}{\varepsilon}\left(\underbrace{S_\alpha-\frac{\sigma}{1+\sigma}S_\alpha C_\alpha-\frac{1}{1+\sigma}C_\alpha}_{=\,O{\left(1\right)}\text{ uniformly as }\varepsilon\rightarrow 0^+}\right)}_{=\,O{\left(\dfrac{1}{\varepsilon}\right)}\text{  uniformly as }\varepsilon\rightarrow 0^+}$$ as well as $$\underbrace{\dfrac{\de C_\alpha}{\de t_\alpha}}_{\substack{=\,O{\left(1\right)}\text{  uniformly as }\varepsilon\rightarrow 0^+,\\ \text{when }t_\alpha\,=\,O{\left(1\right)}\text{  uniformly as }\varepsilon\rightarrow 0^+}}=\underbrace{\dots}_{=\,O{\left(\dfrac{1}{\varepsilon}\right)}\text{  uniformly as }\varepsilon\rightarrow 0^+},$$ therefore $$\underbrace{S_\alpha-\frac{\sigma}{1+\sigma}S_\alpha C_\alpha-\frac{1}{1+\sigma}C_\alpha}_{\substack{=\,O{\left(\varepsilon\right)}\text{  uniformly as }\varepsilon\rightarrow 0^+,\\ \text{when }t_\alpha\,=\,O{\left(1\right)}\text{  uniformly as }\varepsilon\rightarrow 0^+}}.$$ Thus, it follows that
\begin{equation}
\label{deCalphat2approx}
\begin{split}
\varepsilon\frac{\de C_\alpha}{\de t_\alpha}=S_\alpha-\frac{\sigma}{1+\sigma}S_\alpha C_\alpha-\frac{1}{1+\sigma}C_\alpha\approx 0,&\\ &\hspace{-3cm}\text{when }\exists B>0\text{ independent of }\varepsilon,\text{ such that: }\left|t_\alpha\right|\leq B,
\end{split}
\end{equation}
i.e., $$C_\alpha\approx\frac{\left(1+\sigma\right)S_\alpha}{1+\sigma S_\alpha},\text{ when }\exists B>0\text{ independent of }\varepsilon,\text{ such that: }\left|t_\alpha\right|\leq B.$$ 
If we insert the above approximate equality in  (\ref{SCeqdmlsst2;a}), then the later becomes an approximate separable non linear differential  equation, i.e., 
\begin{equation}
\label{deSalphat2approx}
\frac{\de S_\alpha}{\de t_\alpha}\approx -\frac{S_\alpha}{\left(1+\rho\right)\left(1+\sigma S_\alpha\right)},\text {when }\exists B>0\text{ independent of }\varepsilon,\text{ such that: }\left|t_\alpha\right|\leq B,
\end{equation}
the solution of which is
\begin{equation*}
S_\alpha\approx \frac{1}{\sigma}W{\left(\sigma\ell
\exp{\left\{\sigma \ell-\frac{1}{1+\rho}\, t_\alpha\right\}}\right)},\text{ when }\exists B>0:\text{ }\left|t_\alpha\right|\leq B,
\end{equation*}
where $W$ is the aforementioned Lambert function, and $\ell\geq 0$ is a constant that remains to be determined.\\ Therefore, the outer approximation, $\left({\left[S\right]}_{{\rm out}},{\left[C\right]}_{{\rm out}}\right)$, of the solution $\left(\left[S\right],\left[C\right]\right)$ of (\ref{SC_Cauchy}), i.e., the approximation for those $t$ for which it holds that $$\text{\fbox{$\exists B>0$ independent of $\varepsilon$, such that: $\left|t\right|\leq B t_2$}},$$ is 
\begin{equation*}
\begin{split}
{\left[S\right]}_{{\rm out}}&=K_{\!M} W{\left(\frac{\ell A_1}{K_{\!M}}\exp{\left\{\frac{1}{K_{\!M}}\left(\ell A_1-k_2 A_2 t\right)\right\}}\right)}\\
{\left[C\right]}_{{\rm out}}&=\frac{A_2 W{\left(\frac{\ell A_1}{K_{\!M}}\exp{\left\{\frac{1}{K_{\!M}}\left(\ell A_1-k_2 A_2 t\right)\right\}}\right)}}{1+W{\left(\frac{\ell A_1}{K_{\!M}}\exp{\left\{\frac{1}{K_{\!M}}\left(\ell A_1-k_2 A_2 t\right)\right\}}\right)}}.
\end{split}
\end{equation*}
\end{itemize}
We can now utilise the matching technique in order to find a uniform approximation of the solution
$\left(\left[S\right],\left[C\right]\right)$ of (\ref{SC_Cauchy}) from the individual approximations $\left({\left[S\right]}_{{\rm in}},{\left[C\right]}_{{\rm in}}\right)$ and $\left({\left[S\right]}_{{\rm out}},{\left[C\right]}_{{\rm out}}\right)$. First, choosing a time scale between $t_1$ and $t_2$, e.g. 
$$\frac{t_1}{\varepsilon^{\frac{1}{2}}}\in\left(t_1,t_2\right),$$  
we find easily that the common limit resulting from the matching condition of the two individual solutions should be  $$L\coloneqq\left({\left[S\right]}_0,\frac{A_2{\left[S\right]}_0}{K_{\!M}+{\left[S\right]}_0}\right).$$ Therefore $$\ell=\frac{{\left[S\right]}_0}{A_1}$$ and thus a uniform approximation $\left({\left[S\right]}_{{\rm un}},{\left[C\right]}_{{\rm un}}\right)$, of $\left(\left[S\right],\left[C\right]\right)$ is $$\left({\left[S\right]}_{{\rm un}},{\left[C\right]}_{{\rm un}}\right)=\left({\left[S\right]}_{{\rm in}},{\left[C\right]}_{{\rm in}}\right)+\left({\left[S\right]}_{{\rm out}},{\left[C\right]}_{{\rm out}}\right)-L,$$ i.e., in more detail, 
\begin{equation}
\label{SuCu}
\begin{split}
{\left[S\right]}_{{\rm un}}&=K_{\!M} W{\left(\frac{{\left[S\right]}_0}{K_{\!M}}\exp{\left\{\frac{1}{K_{\!M}}\left({\left[S\right]}_0-k_2 A_2 t\right)\right\}}\right)},\\
{\left[C\right]}_{{\rm un}}&=\frac{A_2 W{\left(\frac{{\left[S\right]}_0}{K_{\!M}}\exp{\left\{\frac{1}{K_{\!M}}\left({\left[S\right]}_0-k_2 A_2 t\right)\right\}}\right)}}{1+W{\left(\frac{{\left[S\right]}_0}{K_{\!M}}\exp{\left\{\frac{1}{K_{\!M}}\left({\left[S\right]}_0-k_2 A_2 t\right)\right\}}\right)}}+\left({\left[C\right]}_0-\frac{A_2{\left[S\right]}_0}{K_{\!M}+{\left[S\right]}_0}\right)\exp{\left\{-k_1\left(K_{\!M}+{\left[S\right]}_0\right)t  \right\}}.
\end{split}
\end{equation}

\subsection{Total substrate approach} 

Since $\left(\left[S\right],\left[C\right]\right)\in\Omega_1$, then $\left[T\right]\leq A_1$, where $\left[T\right]$ is as in (\ref{TclneqSC}).
Thus, we introduce, as usually, the dimensionless dependent variable
 $$T_\alpha{\left(t_\alpha\right)}\coloneqq\frac{1}{A_1}{\left(\left[S\right]+\left[C\right]\right)}{\left(t_\alpha\right)}={\left(S_\alpha+\varepsilon C_\alpha\right)}{\left(t_\alpha\right)}.$$ 
 It is easily verified that (\ref{SCeqdmlsst1}) and (\ref{SCeqdmlsst2}) will take the following forms:
\begin{itemize}
\item If \fbox{$t_\alpha=\dfrac{t}{t_1}$}, then 
\begin{subequations}
\label{TCeqdmlsst1}
\begin{align}
\dfrac{\de T_\alpha}{\de t_\alpha}&=-\frac{\varepsilon}{\left(1+\rho\right)\left(1+\sigma\right)}C_\alpha\label{TCeqdmlsst1;a},\\
\dfrac{\de C_\alpha}{\de t_\alpha}&=\frac{\varepsilon\sigma}{1+\sigma}{C_\alpha}^2-\left(\frac{1}{1+\sigma}+\varepsilon\right)C_\alpha-\frac{\sigma}{1+\sigma}C_\alpha T_\alpha+T_\alpha\label{TCeqdmlsst1;b}. 
\end{align}
\end{subequations}
\item If \fbox{$t_\alpha=\dfrac{t}{t_2}$}, then 
\begin{subequations}
\label{TCeqdmlsst2}
\begin{align}
\dfrac{\de T_\alpha}{\de t_\alpha}&=-\frac{1}{\left(1+\rho\right)\left(1+\sigma\right)}C_\alpha\label{TCeqdmlsst2;a},\\
\dfrac{\de C_\alpha}{\de t_\alpha}&=\frac{\sigma}{1+\sigma}{C_\alpha}^2-\left(\frac{1}{\varepsilon\left(1+\sigma\right)}+1\right)C_\alpha-\frac{\sigma}{\varepsilon\left(1+\sigma\right)}C_\alpha T_\alpha+\frac{1}{\varepsilon}T_\alpha\label{TCeqdmlsst2;b}. 
\end{align}
\end{subequations}
\end{itemize}
So we have the following scaled problem:
\[\label{TCdmlss_Cauchy}\tag{$TC\alpha_s$}
\begin{minipage}{.8\textwidth}
\centering 
\textit{Given ${\left[S\right]}_0,{\left[C\right]}_0,{\left[P\right]}_0\,\geq 0$ and $\varepsilon>0$, we seek an interval $\mathcal{I}\subseteq\mathbb{R}$ with $0\in\mathcal{I}$, and a function $\left(T_\alpha,C_\alpha\right)\colon\mathcal{I}\rightarrow {\left[0,1\right]}^2$, such that $\left(T_\alpha,C_\alpha\right)$ satisfies both \eqref{TCeqdmlsst1} if $t_\alpha=\dfrac{t}{t_1}$, or \eqref{TCeqdmlsst2} if $t_\alpha=\dfrac{t}{t_2}$ in ${\left(\mathcal{I}\setminus\left\{0\right\}\right)}^\circ$, and $\left(T_\alpha,C_\alpha\right)=\left(\dfrac{{\left[S\right]}_0}{A_1}+\dfrac{{\left[C\right]}_0}{\varepsilon A_1},\dfrac{{\left[C\right]}_0}{\varepsilon A_1}\right)$, for $t_\alpha=0$.}
\end{minipage}
\]
Working as with problem (\ref{SCdmlss_Cauchy}), we conclude for problem (\ref{TCdmlss_Cauchy}) now, the following: 
\begin{itemize}
\item (\ref{TCeqdmlsst1;a}) gives that $$
\frac{\de T_\alpha}{\de t_\alpha}\approx 0,\text{ when }\exists B>0\text{ independent of }\varepsilon,\text{ such that: }\left|t_\alpha\right|\leq B,$$
and due to the initial condition of (\ref{TCdmlss_Cauchy}) we have that $$T_\alpha\approx \frac{{\left[S\right]}_0}{A_1}+\frac{{\left[C\right]}_0}{\varepsilon A_1},\text{ when }\exists B>0:\text{ }\left|t_\alpha\right|\leq B.$$ Inserting the above approximate equality into (\ref{TCeqdmlsst1;b}), which in turn takes the following approximate form
\begin{equation*}
\begin{split}
\frac{\de C_\alpha}{\de t_\alpha}\approx-\left(\frac{\sigma}{1+\sigma}T_\alpha+\frac{1}{1+\sigma}\right)C_\alpha+T_\alpha,&\\ &\hspace{-3cm}\text{when }\exists B>0\text{ independent of }\varepsilon,\text{ such that: }\left|t_\alpha\right|\leq B,
\end{split}
\end{equation*}
then the later becomes an approximate linear differential  equation, the solution of which is $$C_\alpha\approx \frac{\left(1+\sigma\right)T_\alpha}{1+\sigma T_\alpha}+\left(\dfrac{{\left[C\right]}_0}{\varepsilon A_1}-\frac{\left(1+\sigma\right)T_\alpha}{1+\sigma T_\alpha}\right)\exp{\left\{-\frac{1+\sigma T_\alpha}{1+\sigma}\, t_\alpha\right\}},\text{ when }\exists B>0:\text{ }\left|t_\alpha\right|\leq B,$$
given the initial condition of (\ref{TCdmlss_Cauchy}).\\ Therefore, the inner approximation, $\left({\left[T\right]}_{{\rm in}},{\left[C\right]}_{{\rm in}}\right)$, of $\left(\left[T\right],\left[C\right]\right)$ is $$\left({\left[T\right]}_{{\rm in}},{\left[C\right]}_{{\rm in}}\right)=\left({\left[T\right]}_0,\frac{A_2{\left[T\right]}_0}{K_{\!M}+{\left[T\right]}_0}+\left({\left[C\right]}_0-\frac{A_2{\left[T\right]}_0}{K_{\!M}+{\left[T\right]}_0}\right)\exp{\left\{-k_1\left(K_{\!M}+{\left[T\right]}_0\right)t  \right\}}\right),$$ where $${\left[T\right]}_0\coloneqq {\left[S\right]}_0+{\left[C\right]}_0.$$ 
\item From (\ref{TCeqdmlsst2;b}) we get that 
\begin{equation*}
\begin{split}
\varepsilon\left(\frac{\de C_\alpha}{\de t_\alpha}-\frac{\sigma}{1+\sigma}{C_\alpha}^2+C_\alpha\right)=-\frac{1}{1+\sigma}C_\alpha-\frac{\sigma}{1+\sigma}C_\alpha T_\alpha+T_\alpha\approx 0,&\\ &\hspace{-5cm}\text{when }\exists B>0\text{ independent of }\varepsilon,\text{ such that: }\left|t_\alpha\right|\leq B,
\end{split}
\end{equation*}
i.e.,
$$C_\alpha\approx\frac{\left(1+\sigma\right)T_\alpha}{1+\sigma T_\alpha},\text{ when }\exists B>0\text{ independent of }\varepsilon,\text{ such that: }\left|t_\alpha\right|\leq B.$$ 
If we insert the above approximate equality into (\ref{TCeqdmlsst2;a}), then the later becomes an approximate separable non linear differential  equation, namely $$\frac{\de T_\alpha}{\de t_\alpha}\approx -\frac{T_\alpha}{\left(1+\rho\right)\left(1+\sigma T_\alpha\right)},\text { when }\exists B>0\text{ independent of }\varepsilon,\text{ such that: }\left|t_\alpha\right|\leq B,$$ the solution of which is $$T_\alpha\approx \frac{1}{\sigma}W\left(\sigma\ell
\exp{\left\{\sigma \ell-\frac{1}{1+\rho}\, t_\alpha\right\}}\right),\text{ when }\exists B>0:\text{ }\left|t_\alpha\right|\leq B,$$ where $\ell\geq 0$ is a constant that remains to be determined. Therefore, the outer approximation $\left({\left[T\right]}_{{\rm out}},{\left[C\right]}_{{\rm out}}\right)$ of $\left(\left[T\right],\left[C\right]\right)$, is 
\begin{equation*}
\begin{split}
{\left[T\right]}_{{\rm out}}&=K_{\!M} W{\left(\frac{\ell A_1}{K_{\!M}}\exp{\left\{\frac{1}{K_{\!M}}\left(\ell A_1-k_2 A_2 t\right)\right\}}\right)}\\
{\left[C\right]}_{{\rm out}}&=\frac{A_2 W{\left(\frac{\ell A_1}{K_{\!M}}\exp{\left\{\frac{1}{K_{\!M}}\left(\ell A_1-k_2 A_2 t\right)\right\}}\right)}}{1+W{\left(\frac{\ell A_1}{K_{\!M}}\exp{\left\{\frac{1}{K_{\!M}}\left(\ell A_1-k_2 A_2 t\right)\right\}}\right)}}.
\end{split}
\end{equation*}
\end{itemize}

Finally, with a similar reasoning as for the uniform approximation of the solution of (\ref{SC_Cauchy}), we have that $$\ell=\frac{{\left[T\right]}_0}{A_1}$$ and also that the uniform approximation, $\left({\left[T\right]}_{{\rm un}},{\left[C\right]}_{{\rm un}}\right)$, of $\left(\left[T\right],\left[C\right]\right)$ is 
\begin{equation}
\label{TuCu}
\begin{split}
{\left[T\right]}_{{\rm un}}&=K_{\!M} W{\left(\frac{{\left[T\right]}_0}{K_{\!M}}\exp{\left\{\frac{1}{K_{\!M}}\left({\left[T\right]}_0-k_2 A_2 t\right)\right\}}\right)}\\
{\left[C\right]}_{{\rm un}}&=\frac{A_2 W{\left(\frac{{\left[T\right]}_0}{K_{\!M}}\exp{\left\{\frac{1}{K_{\!M}}\left({\left[T\right]}_0-k_2 A_2 t\right)\right\}}\right)}}{1+W{\left(\frac{{\left[T\right]}_0}{K_{\!M}}\exp{\left\{\frac{1}{K_{\!M}}\left({\left[T\right]}_0-k_2 A_2 t\right)\right\}}\right)}}+\left({\left[T\right]}_0-\frac{A_2{\left[T\right]}_0}{K_{\!M}+{\left[T\right]}_0}\right)\exp{\left\{-k_1\left(K_{\!M}+{\left[T\right]}_0\right)t  \right\}}.
\end{split}
\end{equation}

\subsection{Conclusions} 

We showed that given (\ref{sQSSA}) there is a $t_1>0$ such that $$\frac{\de \left[S\right]}{\de t}\approx 0,\text{ when }\exists B>0\text{ independent of }\varepsilon,\text{ such that: }\left|t\right|\leq B t_1,$$ which arises directly from (\ref{deSalphat1approx}), as well as that there is a $t_2\gg t_1$ such that $$\frac{\de \left[C\right]}{\de t}\approx 0,\text{ when }\exists B>0\text{ independent of }\varepsilon,\text{ such that: }\left|t\right|\leq B t_2,$$ which in turn results from (\ref{deCalphat2approx}).\\
In fact, due to (\ref{deSalphat2approx}) it holds that $$\frac{\de \left[S\right]}{\de t}\approx -\frac{k_2 A_2\left[S\right]}{K_{\!M}+\left[S\right]},\text{ when }\exists B>0:\text{ }\left|t\right|\leq B t_2.$$ 
Hence, we can conclude that $$\upsilon\approx\begin{cases}
0,&\text{ when }\exists B>0:\text{ }\left|t\right|\leq B t_1,\\
\dfrac{k_2 A_2 \left[S\right]}{K_{\!M}+\left[S\right]},&\text{ when }\exists B>0:\text{ }\left|t\right|\leq B t_2,
\end{cases}$$
where $\upsilon$ is the rate of the chemical reaction with chemical equation (\ref{SEP}). The above approximation for $t$ comparable to $t_2$ is none other than the Michaelis-Menten  approximation for the kinetics of the aforementioned chemical reaction, as already commented in (\ref{MMeq}).

\begin{figure}[!htbp]
\centering
\includegraphics[width=.5\textwidth]{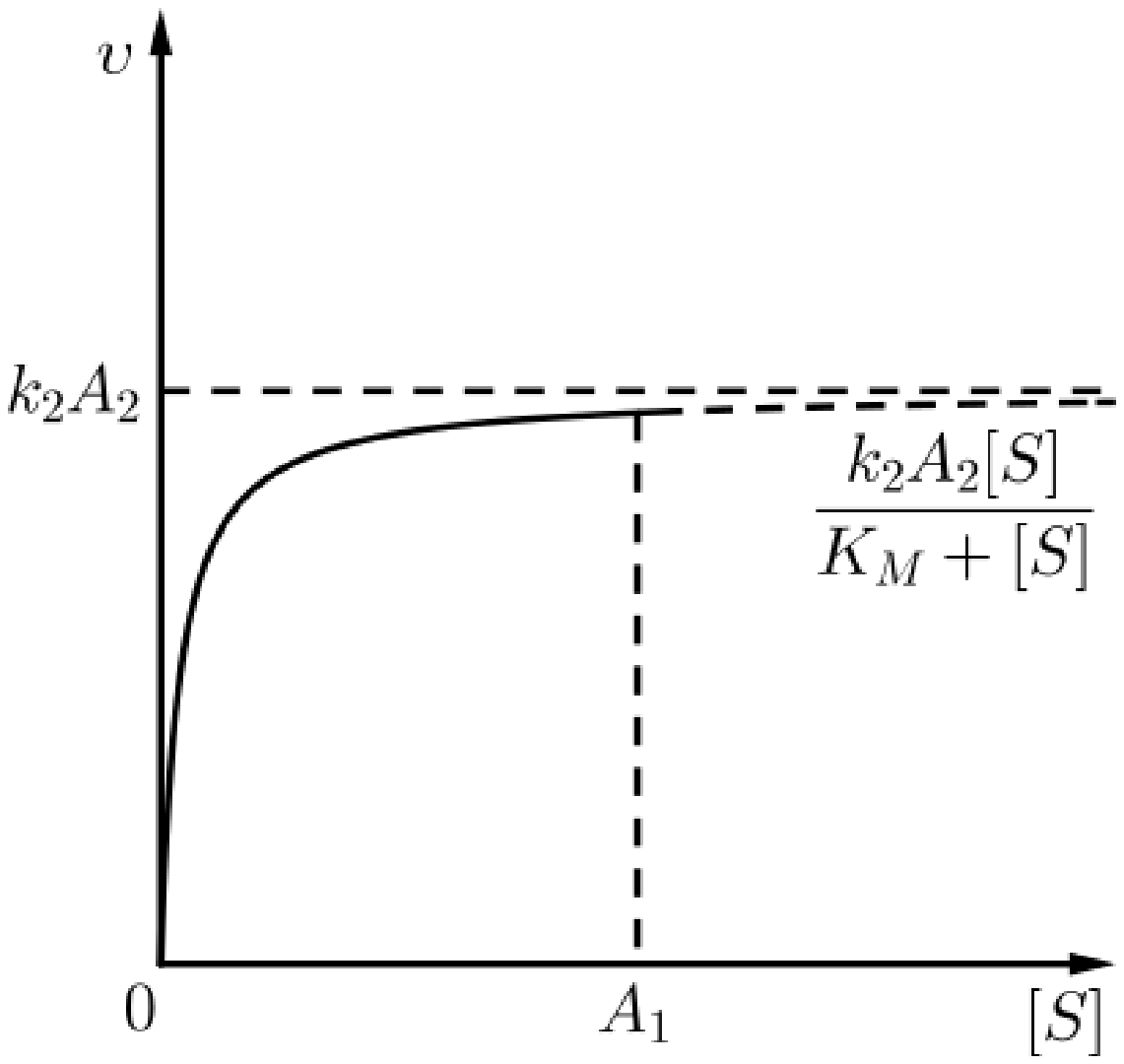}
\caption{An approximation for the kinetics of the chemical reaction (\ref{SEP}) given that (\ref{sQSSA}) holds, for times comparable to $t_2$.}
\label{MMupsilon}
\end{figure}

Furthermore, comparing the approximate solution of the free substrate $\left({\left[S\right]}_{{\rm un}},{\left[C\right]}_{{\rm un}}\right)$ of (\ref{SuCu}) with the approximate solution of the total substrate $\left({\left[T\right]}_{{\rm un}},{\left[C\right]}_{{\rm un}}\right)$ of (\ref{TuCu}), these two should be in agreement. Indeed, it is sufficient to observe that $$\left[T\right]\approx\left[S\right],$$ as $$T_\alpha=S_\alpha+\varepsilon C_\alpha\overset{\text{(\ref{epsln2})}}{\approx}S_\alpha\Rightarrow T\approx S.$$

\subsection{Numerical solution} 

\begin{figure}[!htbp]
\begin{subfigure}{.5\textwidth}
  \centering
\includegraphics[width=1\linewidth]{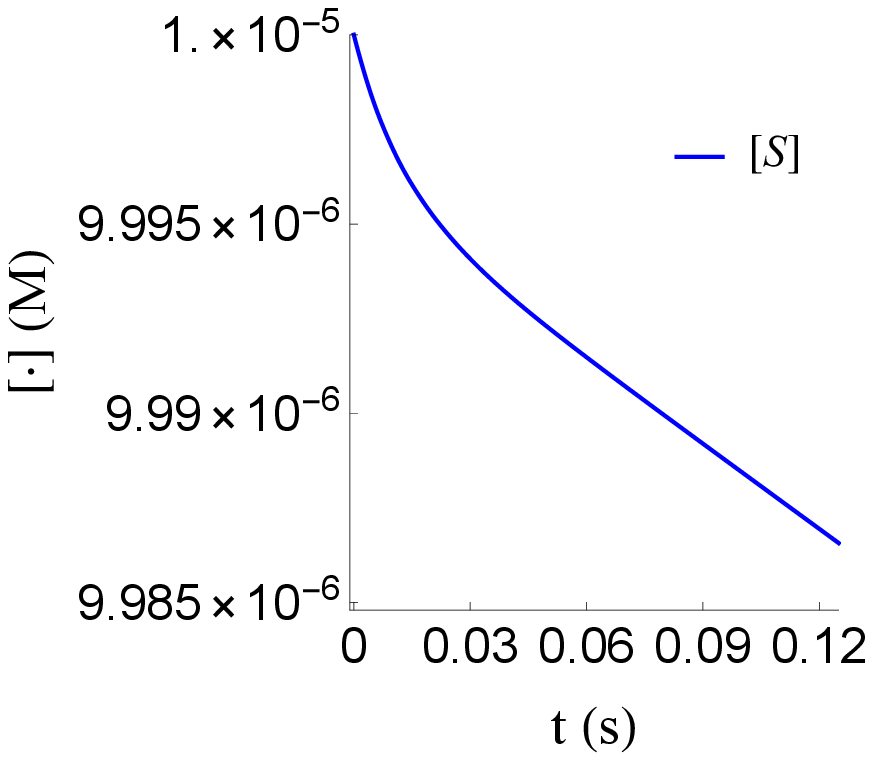}
 \caption{}
  \label{sQSSA-numer-sol.a}
\end{subfigure}
\begin{subfigure}{.5\textwidth}
  \centering
 \includegraphics[width=1\linewidth]{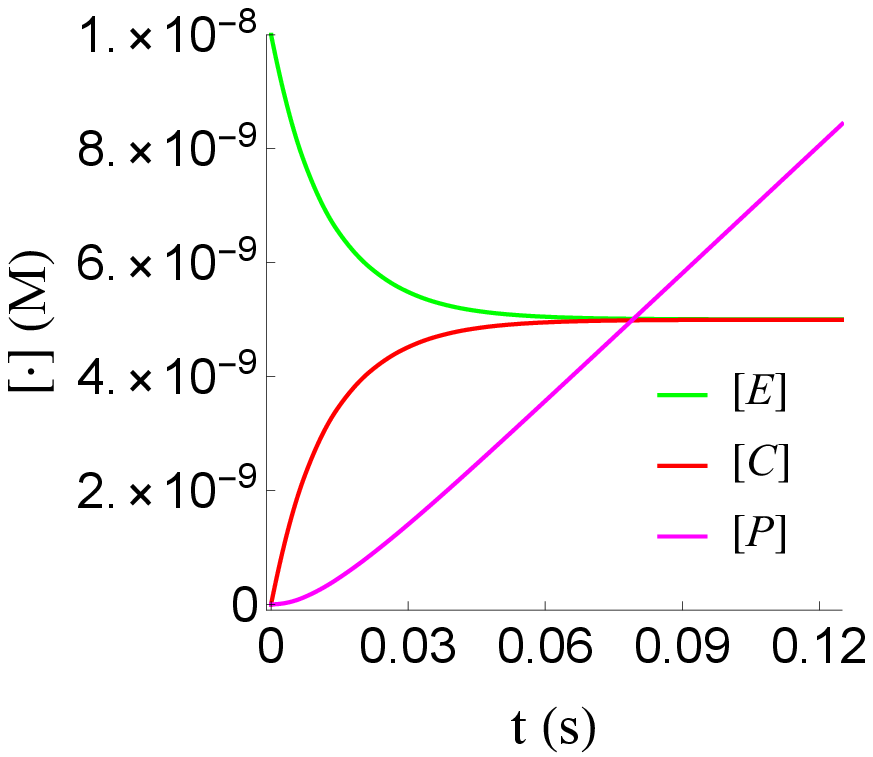}
  \caption{}
  \label{sQSSA-numer-sol.b}
\end{subfigure}
\begin{subfigure}{.5\textwidth}
  \centering
   \includegraphics[width=1\linewidth]{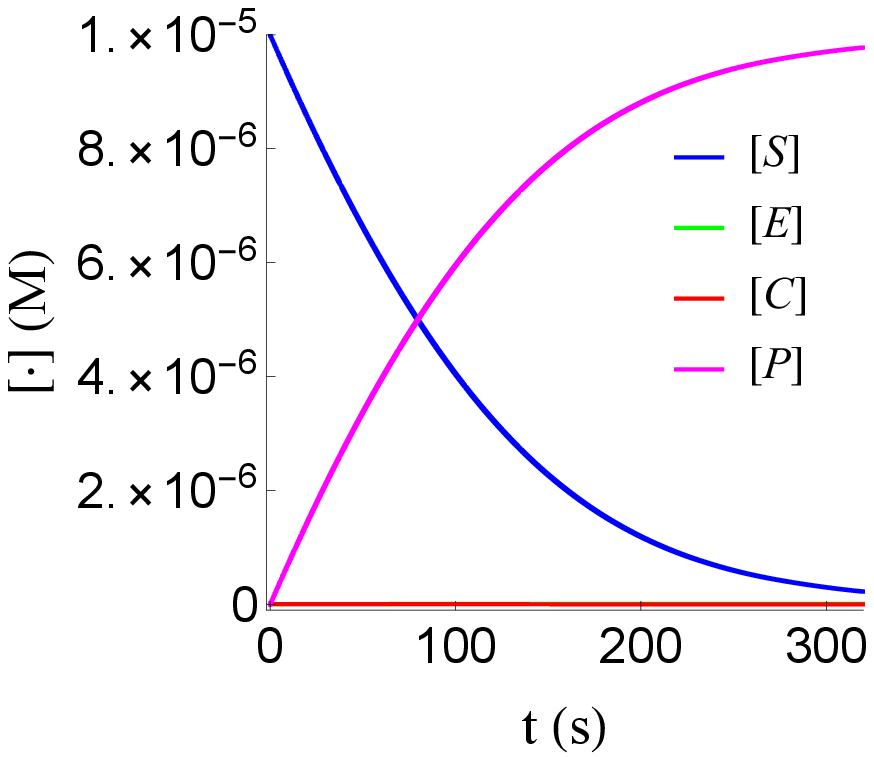}  
  \caption{}
  \label{sQSSA-numer-sol.c}
\end{subfigure}
\begin{subfigure}{.5\textwidth}
  \centering
\includegraphics[width=1\linewidth]{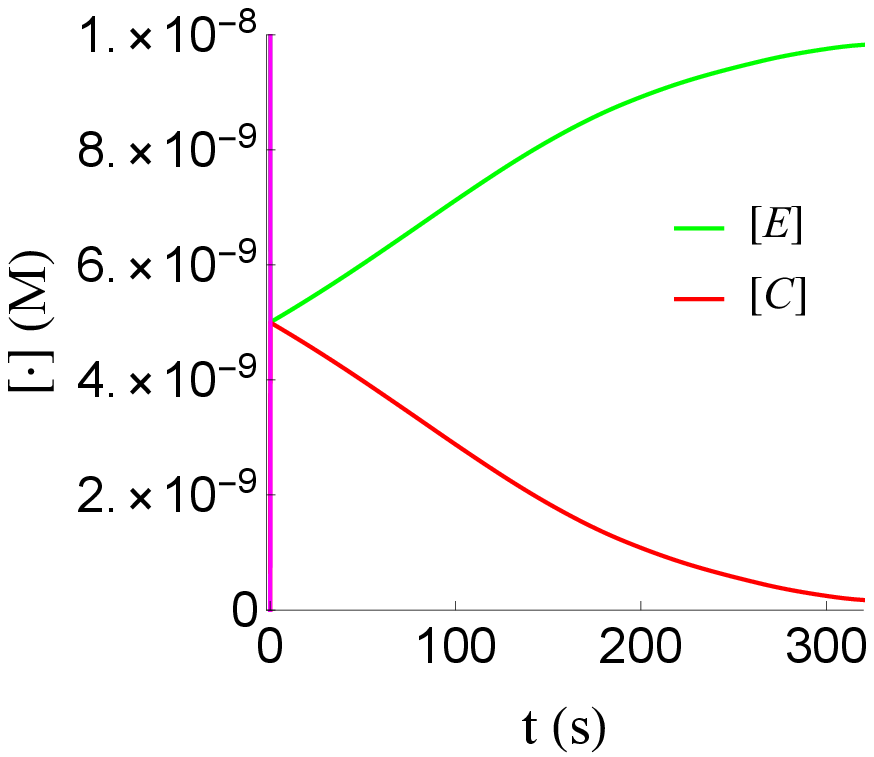}  
  \caption{}
  \label{sQSSA-numer-sol.d}
\end{subfigure}
  \caption{Plot of $\left[S\right]$, $\left[E\right]$, $\left[C\right]$ and $\left[P\right]$ of problem (\ref{SECP_Cauchy}) for non negative times, given that \eqref{sQSSA} holds. We see that $\left[S\right]$ and $\left[C\right]$ are of different order of magnitude, as well as that there are two distinct phases of the evolution of the phenomenon.}
   \label{sQSSA-numer-sol}
\end{figure}

\begin{figure}[!htbp]
\begin{subfigure}{.5\textwidth}
  \centering
\includegraphics[width=1\linewidth]{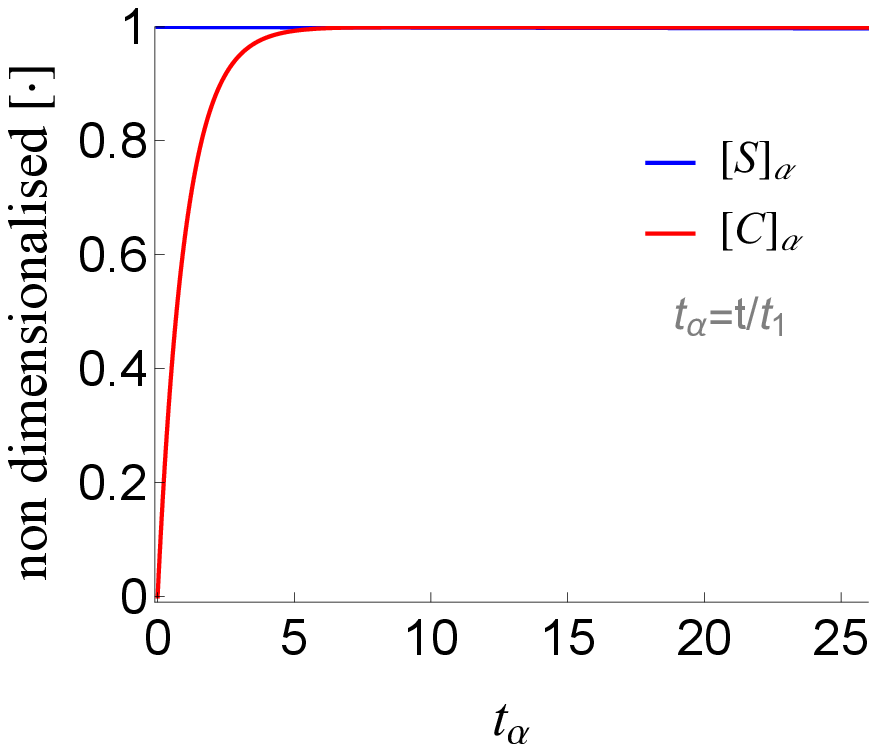}
 \caption{}
  \label{sQSSA-numer-sol-norm.a}
\end{subfigure}
\begin{subfigure}{.5\textwidth}
  \centering
 \includegraphics[width=1\linewidth]{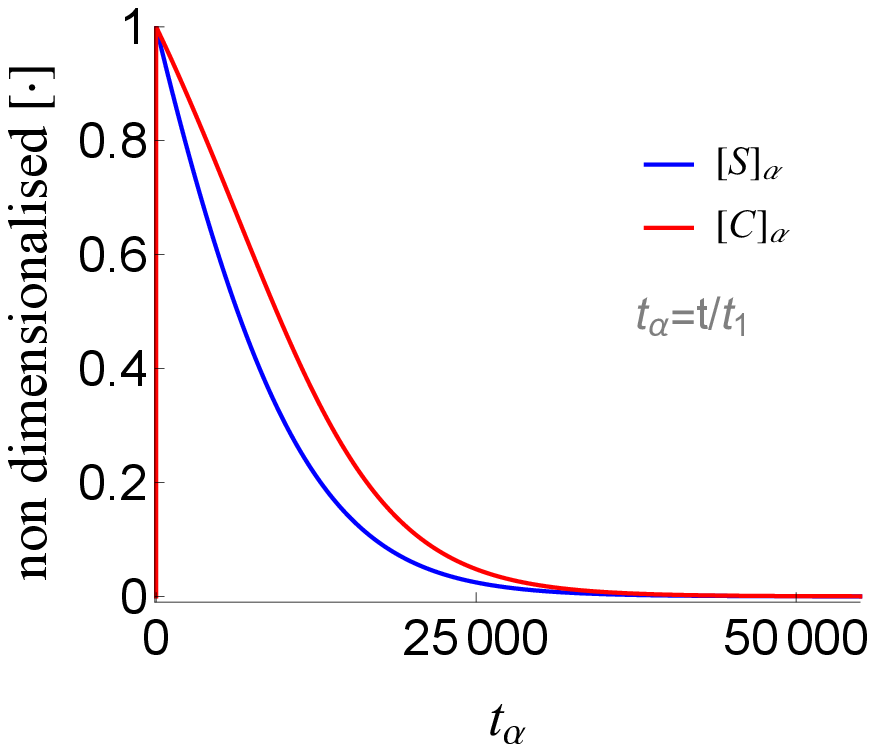}
  \caption{}
  \label{sQSSA-numer-sol-norm.b}
\end{subfigure}
\begin{subfigure}{.5\textwidth}
  \centering
   \includegraphics[width=1\linewidth]{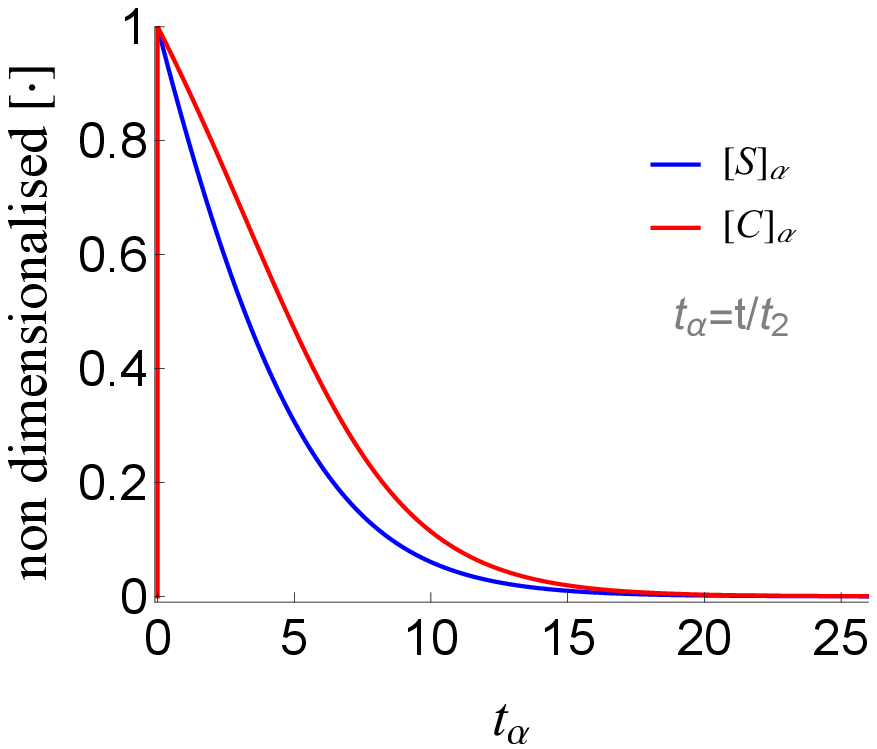}  \caption{}
  \label{sQSSA-numer-sol-norm.c}
\end{subfigure}
\begin{subfigure}{.5\textwidth}
  \centering
\includegraphics[width=1\linewidth]{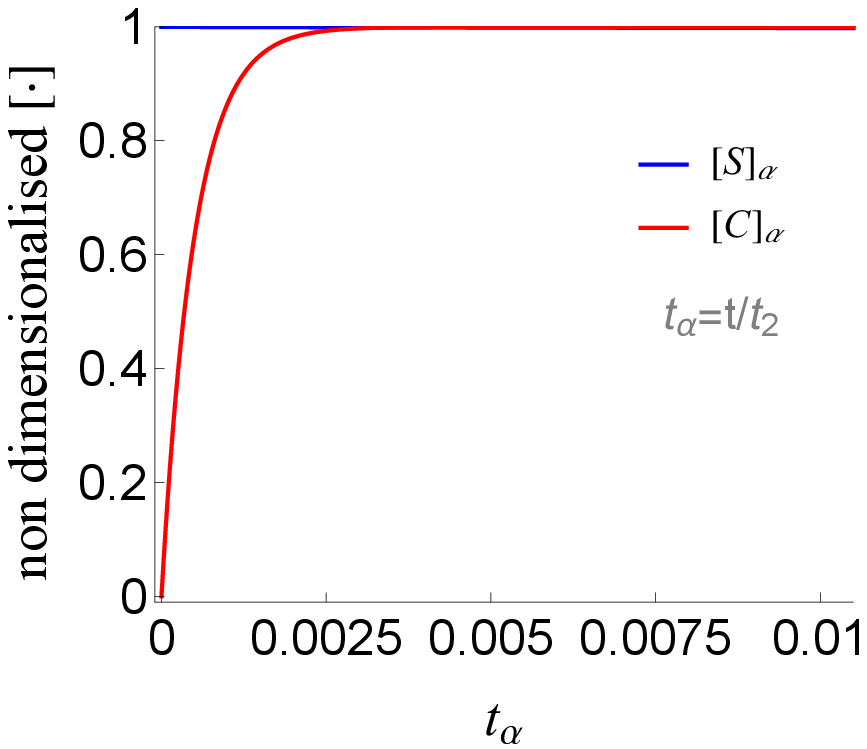}  
  \caption{}
  \label{sQSSA-numer-sol-norm.d}
\end{subfigure}
  \caption{Plots of $S_\alpha$ and $C_\alpha$ of problem (\ref{SCdmlss_Cauchy}) for non negative times, given that \eqref{sQSSA} holds. In (a) and (b) time is measured based on $t_1$, whereas in (c) and (d) based on $t_2$.}
   \label{sQSSA-numer-sol-norm}
\end{figure}

\begin{figure}[!htbp]
\begin{subfigure}{.5\textwidth}
  \centering
\includegraphics[width=1\linewidth]{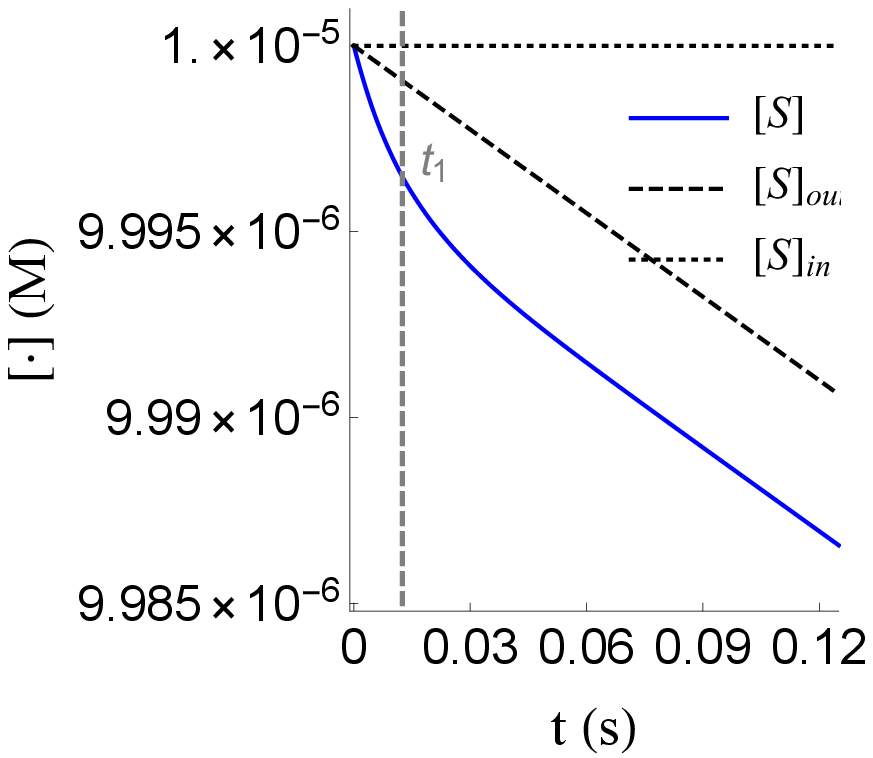}
 \caption{}
  \label{sQSSA-numer-sol-approx.a}
\end{subfigure}
\begin{subfigure}{.5\textwidth}
  \centering
 \includegraphics[width=1\linewidth]{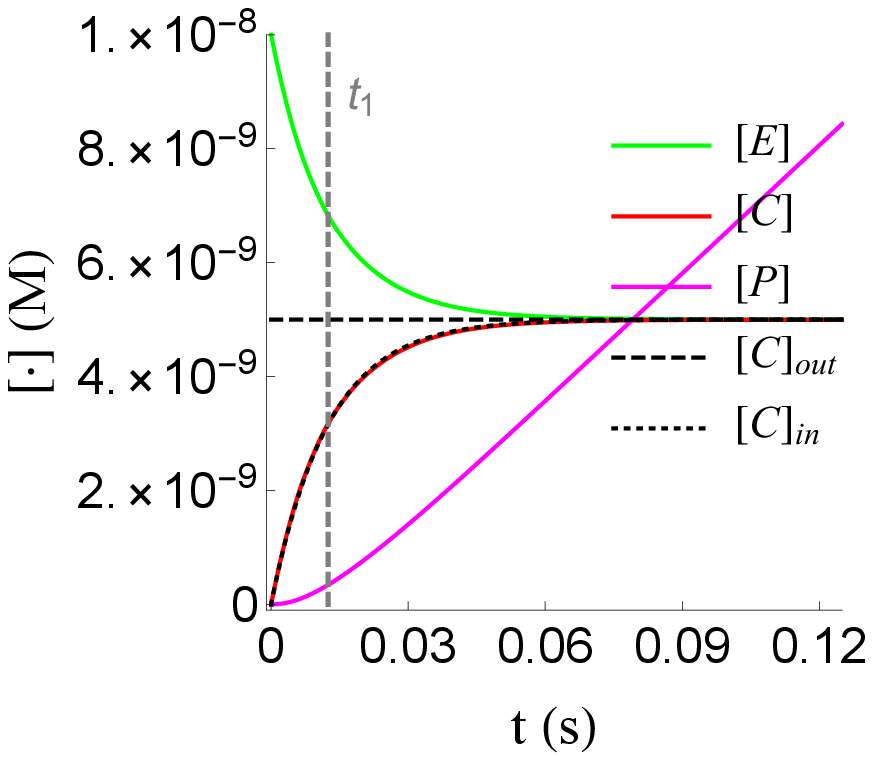}
  \caption{}
  \label{sQSSA-numer-sol-approx.b}
\end{subfigure}
\begin{subfigure}{.5\textwidth}
  \centering
   \includegraphics[width=1\linewidth]{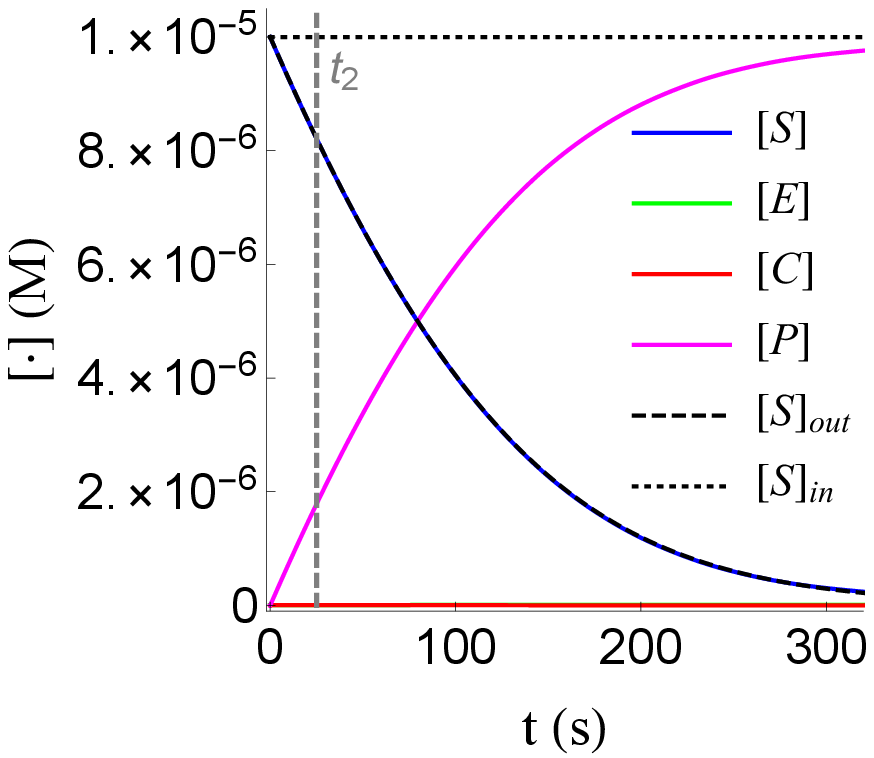}  
  \caption{}
  \label{sQSSA-numer-sol-approx.c}
\end{subfigure}
\begin{subfigure}{.5\textwidth}
  \centering
\includegraphics[width=1\linewidth]{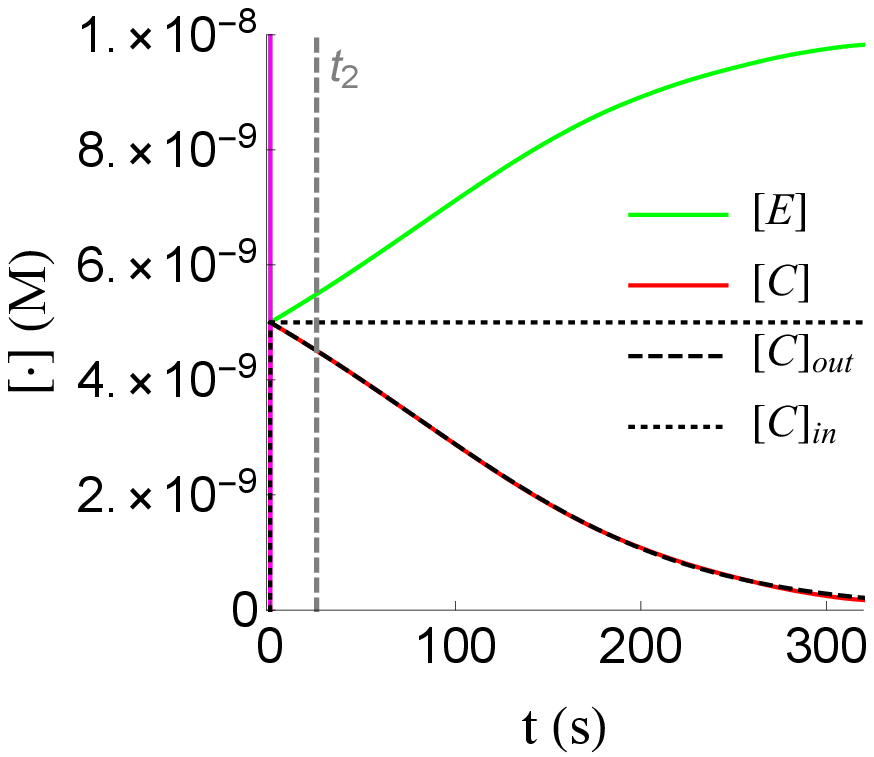}  
  \caption{}
  \label{sQSSA-numer-sol-approx.d}
\end{subfigure}
  \caption{Plots of the inner and outer approximations of ${\left[S\right]}$ and $\left[C\right]$ of problem (\ref{SECP_Cauchy}), for non negative times, given that \eqref{sQSSA} holds.}
   \label{sQSSA-numer-sol-approx}
\end{figure}

We proceed to the numerical solution of the problem, as shown in \hyperref[sQSSA-numer-sol]{Figure \ref*{sQSSA-numer-sol}}, in \hyperref[sQSSA-numer-sol-norm]{Figure \ref*{sQSSA-numer-sol-norm}} and in \hyperref[sQSSA-numer-sol-approx]{Figure \ref*{sQSSA-numer-sol-approx}}, to verify our conclusions. For the numerical values of the constants and the initial conditions we follow the work of Segel in 1988 \cite{segel1988validity}. The values are given in the table below. 

\begin{table}[ht]
\centering 
\begin{tabular}{c c c}
\hline
Parameter & Value & Unit \\ [0.5ex]
\hline\hline          
$k_{-1}$  & $25$ & $s^{-1}$ \\
$k_1$  & $4\cdot 10^6$ & $M^{-1}s^{-1}$ \\
$k_2$  & $15$ & $s^{-1}$ \\
${\left[S\right]}_0$  & $10^{-5}$ & $M$\\
${\left[E\right]}_0$  & $10^{-8}$ & $M$ \\
${\left[C\right]}_0$  & $0$ & $M$ \\
${\left[P\right]}_0$  & $0$ & $M$ \\  [1ex]      
\hline
\end{tabular}
\label{table:2}
\end{table}
\noindent We calculate $$K_{\!M} ={\left[S\right]}_0=A_1\text{ and }A_2={\left[E\right]}_0=10^{-3}{\left[S\right]}_0=10^{-3}A_1,$$ i.e., $$\varepsilon=5\cdot 10^{-4}$$ and $$t_1=1,25\cdot 10^{-2}s\text{ and }t_2=25s.$$

\section{The reverse quasi-steady-state assumption}

The reverse quasi-steady-state assumption is the following:
\begin{equation}
\label{rQSSA}\tag{$rQSSA$}
\fbox{$A_1>0$ and $A_2\gg K_{\!M}+A_1$},
\end{equation} 
or, equivalently, $$A_1>0,\text{ }A_2\gg K_{\!M}\text{ and }A_2\gg A_1,$$ and given that this holds we study problem (\ref{SC_Cauchy}).\\
We notice that the  inequality $A_2\gg K_{\!M}+A_1$ of (\ref{rQSSA}) is the reverse of the one corresponding to (\ref{sQSSA}), hence the name of the first. As for the analysis of (\ref{sQSSA}), here, as well, we consider two approaches for the study of case  (\ref{rQSSA}), the free substrate approach and the total substrate approach. 

\subsection{Free substrate approach}

Using only  (\ref{rQSSA}) we will show that:
\begin{enumerate}
\item Problem  (\ref{SC_Cauchy}), and therefore problem  (\ref{SECP_Cauchy}) as well, has inherently two time scales which we will determine. In particular, except for a short initial time interval where the enzymatic reaction with chemical equation (\ref{SEP}) is evolving with rate $\upsilon$, showing approximately linear behaviour with respect to $\left[S\right]$, and $\upsilon\approx k_1 A_2\left[S\right]$ as in (\ref{upslnssr}), during the rest of the time the enzymatic reaction does not evolve, i.e. $\upsilon\approx 0$.  
\item There is a good uniform approximation in closed form to the solution of (\ref{SC_Cauchy}), and therefore to  (\ref{SECP_Cauchy}) as well,  which we will determine.
\end{enumerate}

As usually, we use the dimensionless dependent variables  $$S_\alpha{\left(t_\alpha\right)}\coloneqq\frac{1}{A_1}\left[S\right]{\left(\frac{t}{t_*}\right)}\text{ and }C_\alpha{\left(t_\alpha\right)}\coloneqq\frac{1}{A_4}\left[C\right]{\left(\frac{t}{t_*}\right)},$$ where we have chosen an arbitrary, for the time being, time scale $t_*>0$ for the scaling.\\ We notice, however, that given (\ref{rQSSA}) it follows from (\ref{bscA12neg}) that $$\fbox{$A_4=A_1$}$$ and $$C_\alpha{\left(t_\alpha\right)}=\frac{1}{A_1}\left[C\right]{\left(\frac{t}{t_*}\right)}.$$ A first conclusion is that the possible change of $\left[S\right]$ is comparable to the corresponding of $\left[C\right]$.\\ We set
\begin{equation}
\label{etavareps}
0\,\overset{\text{(\ref{rQSSA})}}{<}\,\eta\coloneqq\frac{A_1}{A_2}<\frac{K_{\!M}+A_1}{A_2}=\frac{1}{\varepsilon}\,\overset{\text{(\ref{rQSSA})}}{\ll}\,1,
\end{equation}
where $\varepsilon$ is as in (\ref{epsln}), i.e., equivalently 
\begin{equation}
\label{etavareps2}
0<\eta<\frac{1}{\varepsilon}\approx 0,
\end{equation}
and so (\ref{SCeq}) will take the following form
\begin{subequations}
\label{SC2eqdmlss}
\begin{align}
\dfrac{\de S_\alpha}{\de t_\alpha}&=\frac{k_1\left(K_{\!M}+A_1\right)A_2}{A_1}t_*\left(-\frac{\sigma}{1+\sigma}S_\alpha+\frac{\eta\sigma}{1+\sigma}S_\alpha C_\alpha+\frac{\eta\rho}{\left(1+\rho\right)\left(1+\sigma\right)}C_\alpha\right)\label{SC2eqdmlss;a},\\
\dfrac{\de C_\alpha}{\de t_\alpha}&=\frac{k_1\left(K_{\!M}+A_1\right)A_2}{A_1}t_*\left(\frac{\sigma}{1+\sigma}S_\alpha-\frac{\eta\sigma}{1+\sigma}S_\alpha C_\alpha-\frac{\eta}{1+\sigma}C_\alpha\right)\label{SC2eqdmlss;b}, 
\end{align}
\end{subequations}
where $\sigma$ and $\rho$ are as in (\ref{sgmrh}). 

Observing (\ref{SC2eqdmlss}), we define 
\begin{equation}
\label{tStC2}
\frac{A_1}{k_1\left(K_{\!M}+A_1\right)A_2}\eqqcolon t_1\,\overset{\text{(\ref{etavareps})}}{\ll}\, \frac{t_1}{\eta}=\frac{1}{k_1\left(K_{\!M}+A_1\right)}\eqqcolon t_2, 
\end{equation}
to conclude that
\begin{equation*}
t_*=t_1\Rightarrow \frac{k_1\left(K_{\!M}+A_1\right)A_2}{A_1}t_*=1\text{ and }t_*=t_2\Rightarrow \frac{k_1\left(K_{\!M}+A_1\right)A_2}{A_1}t_*=\frac{1}{\eta}, 
\end{equation*}
and so (\ref{SC2eqdmlss}) to get the following forms:
\begin{itemize}
\item If \fbox{$t_\alpha=\dfrac{t}{t_1}$}, then 
\begin{subequations}
\label{SC2eqdmlsst1}
\begin{align}
\dfrac{\de S_\alpha}{\de t_\alpha}&=-\frac{\sigma}{1+\sigma}S_\alpha+\frac{\eta\sigma}{1+\sigma}S_\alpha C_\alpha+\frac{\eta\rho}{\left(1+\rho\right)\left(1+\sigma\right)}C_\alpha\label{SC2eqdmlsst1;a},\\
\dfrac{\de C_\alpha}{\de t_\alpha}&=\frac{\sigma}{1+\sigma}S_\alpha-\frac{\eta\sigma}{1+\sigma}S_\alpha C_\alpha-\frac{\eta}{1+\sigma}C_\alpha\label{SC2eqdmlsst1;b}. 
\end{align}
\end{subequations}
\item If \fbox{$t_\alpha=\dfrac{t}{t_2}$}, then 
\begin{subequations}
\label{SC2eqdmlsst2}
\begin{align}
\dfrac{\de S_\alpha}{\de t_\alpha}&=-\frac{\sigma}{\eta\left(1+\sigma\right)}S_\alpha+\frac{\sigma}{1+\sigma}S_\alpha C_\alpha+\frac{\rho}{\left(1+\rho\right)\left(1+\sigma\right)}C_\alpha\label{SC2eqdmlsst2;a},\\
\dfrac{\de C_\alpha}{\de t_\alpha}&=\frac{\sigma}{\eta\left(1+\sigma\right)}S_\alpha-\frac{\sigma}{1+\sigma}S_\alpha C_\alpha-\frac{1}{1+\sigma}C_\alpha\label{SC2eqdmlsst2;b}. 
\end{align}
\end{subequations}
\end{itemize}

Setting $$\Omega\coloneqq\left\{\left(s,c\right)\in{\left[0,1\right]}^2\,\big|\,s+c\leq 1\right\},$$ the scaled version of (\ref{SC_Cauchy}) will be as follows:
\[\label{SC2dmlss_Cauchy}\tag{$SC\alpha_r$}
\begin{minipage}{.8\textwidth}
\centering 
\textit{Given ${\left[S\right]}_0,{\left[C\right]}_0,{\left[P\right]}_0\,\geq 0$ and $\eta>0$, we seek an interval $\mathcal{I}\subseteq\mathbb{R}$ with $0\in\mathcal{I}$, and a function $\left(S_\alpha,C_\alpha\right)\colon\mathcal{I}\rightarrow \Omega$, such that $\left(S_\alpha,C_\alpha\right)$  satisfies both \eqref{SC2eqdmlsst1} if $t_\alpha=\dfrac{t}{t_1}$, or \eqref{SC2eqdmlsst2} if $t_\alpha=\dfrac{t}{t_2}$ in ${\left(\mathcal{I}\setminus\left\{0\right\}\right)}^\circ$, and $\left(S_\alpha,C_\alpha\right)=\left(\dfrac{{\left[S\right]}_0}{A_1},\dfrac{{\left[C\right]}_0}{A_1}\right)$, for $t_\alpha=0$.}
\end{minipage}
\]
\begin{figure}[!htbp]
\centering
\includegraphics[width=.5\textwidth]{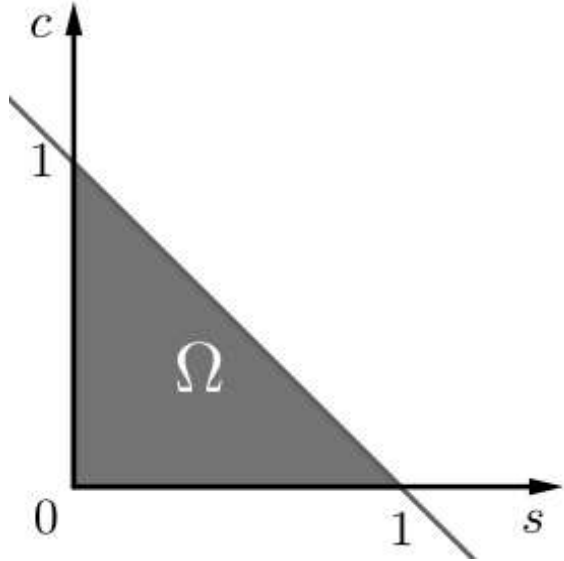}
\caption{The invariant set $\Omega$ of problem \eqref{SC2dmlss_Cauchy}.}
\label{Omega}
\end{figure}

We study each of two versions of (\ref{SC2dmlss_Cauchy}) separately: 
\begin{itemize}
\item From (\ref{SC2eqdmlsst1;a}), which due to (\ref{etavareps2}) takes the approximate linear form  
\begin{equation}
\label{deSalpha2t1approx}
\frac{\de S_\alpha}{\de t_\alpha}\approx -\frac{\sigma}{1+\sigma}S_\alpha,\text{ when }\exists B>0\text{ independent of }\eta,\text{ such that: }\left|t_\alpha\right|\leq B,
\end{equation}
we get, due to the initial condition of (\ref{SC2dmlss_Cauchy}), that $$S_\alpha\approx\frac{{\left[S\right]}_0}{A_1}
\exp{\left\{-\frac{\sigma}{1+\sigma}\, t_\alpha\right\}},\text{ when }\exists B>0:\text{ }\left|t_\alpha\right|\leq B.$$ If we insert the above approximate equality in (\ref{SC2eqdmlsst1;b}), which will now have the approximate form  $$\dfrac{\de C_\alpha}{\de t_\alpha}\approx\frac{\sigma}{1+\sigma}S_\alpha,\text{ when }\exists B>0:\text{ }\left|t_\alpha\right|\leq B,$$ then we  get that $$C_\alpha\approx\frac{{\left[C\right]}_0}{A_1}+\frac{{\left[S\right]}_0}{A_1}\left(1-
\exp{\left\{-\frac{\sigma}{1+\sigma}\, t_\alpha\right\}}\right),\text{ when }\exists B>0:\text{ }\left|t_\alpha\right|\leq B,$$ given the initial condition of (\ref{SC2dmlss_Cauchy}). Therefore, the inner approximation, $\left({\left[S\right]}_{{\rm in}},{\left[C\right]}_{{\rm in}}\right)$, of $\left(\left[S\right],\left[C\right]\right)$ is $$\left({\left[S\right]}_{{\rm in}},{\left[C\right]}_{{\rm in}}\right)=\left({\left[S\right]}_0\exp{\left\{-k_1 A_2 t\right\}},{\left[C\right]}_0+{\left[S\right]}_0\left(1-\exp{\left\{-k_1 A_2 t\right\}}\right)\right).$$ 
\item From (\ref{SC2eqdmlsst2;a}) we have that
\begin{equation}
\label{deSalphat2approxeta}
\begin{split}
\eta\left(\frac{\de S_\alpha}{\de t_\alpha}-\frac{\sigma}{1+\sigma}S_\alpha C_\alpha -\frac{\rho}{\left(1+\rho\right)\left(1+\sigma\right)}C_\alpha\right)=-\frac{\sigma}{1+\sigma}S_\alpha \approx 0,&\\ &\hspace{-6cm}\text{when }\exists B>0\text{ independent of }\eta,\text{ such that: }\left|t_\alpha\right|\leq B,
\end{split}
\end{equation}
i.e., $$S_\alpha \approx 0,\text{ when }\exists B>0:\text{ }\left|t_\alpha\right|\leq B,$$ thus $$\frac{\de S_\alpha}{\de t_\alpha} \approx 0,\text{ when }\exists B>0:\text{ }\left|t_\alpha\right|\leq B.$$
If we insert the above approximate equality in the sum of (\ref{SC2eqdmlsst2;a}) and (\ref{SC2eqdmlsst2;b}), then  the following approximate linear differential  equation arises
 $$\frac{\de C_\alpha}{\de t_\alpha}\approx -\frac{1}{\left(1+\rho\right)\left(1+\sigma\right)}C_\alpha,\text { when }\exists B>0\text{ independent of }\eta,\text{ such that: }\left|t_\alpha\right|\leq B,$$ the solution of which is 
 $$C_\alpha\approx \ell\exp{\left\{-\frac{1}{(1+\rho)(1+\sigma)}\, t_\alpha\right\}},\text{ when }\exists B>0:\text{ }\left|t_\alpha\right|\leq B,$$ where $\ell\geq 0$ a constant that remains to be determined. Therefore, the external approximation, $\left({\left[S\right]}_{{\rm out}},{\left[C\right]}_{{\rm out}}\right)$, of $\left(\left[S\right],\left[C\right]\right)$ is  $$\left({\left[S\right]}_{{\rm out}},{\left[C\right]}_{{\rm out}}\right)=\left(0,\ell A_1\exp{\left\{-k_2 t\right\}}\right).$$ 
\end{itemize}
Finally, as usually, we find that $$\ell=\frac{{\left[S\right]}_0}{A_1}+\frac{{\left[C\right]}_0}{A_1},$$ as well as that the uniform approximation, $\left({\left[S\right]}_{{\rm un}},{\left[C\right]}_{{\rm un}}\right)$, of $\left(\left[S\right],\left[C\right]\right)$ is 
\begin{equation}
\label{SuCur}
\left({\left[S\right]}_{{\rm un}},{\left[C\right]}_{{\rm un}}\right)=\left({\left[S\right]}_0\exp{\left\{-k_1 A_2 t\right\}},{\left[C\right]}_0\exp{\left\{-k_2 t\right\}}+{\left[S\right]}_0\left(\exp{\left\{-k_2 t\right\}}-\exp{\left\{-k_1 A_2 t\right\}}\right)\right).
\end{equation}

\subsection{Total substrate approach} 

Since $\left(\left[S\right],\left[C\right]\right)\in\Omega_1$, we have that $\left[T\right]\leq A_1$. Thus, we introduce, as usually, the dimensionless dependent variable $$T_\alpha{\left(t_\alpha\right)}\coloneqq\frac{1}{A_1}{\left(\left[S\right]+\left[C\right]\right)}{\left(t_\alpha\right)}={\left(S_\alpha+C_\alpha\right)}{\left(t_\alpha\right)}.$$ It is easily verified that (\ref{SC2eqdmlsst1}) and (\ref{SC2eqdmlsst2}) will get the following forms:
\begin{itemize}
\item If \fbox{$t_\alpha=\dfrac{t}{t_1}$}, then 
\begin{subequations}
\label{SCreqdmlsst1}
\begin{align}
\dfrac{\de T_\alpha}{\de t_\alpha}&=-\frac{\eta}{\left(1+\rho\right)\left(1+\sigma\right)}C_\alpha\label{SCreqdmlsst1;a},\\
\dfrac{\de C_\alpha}{\de t_\alpha}&=\frac{\eta\sigma}{1+\sigma}{C_\alpha}^2-\frac{\eta+\sigma}{1+\sigma}C_\alpha-\frac{\eta\sigma}{1+\sigma}C_\alpha T_\alpha+\frac{\sigma}{1+\sigma}T_\alpha\label{SCreqdmlsst1;b}. 
\end{align}
\end{subequations}
\item If \fbox{$t_\alpha=\dfrac{t}{t_2}$}, then 
\begin{subequations}
\label{SCreqdmlsst2}
\begin{align}
\dfrac{\de T_\alpha}{\de t_\alpha}&=-\frac{1}{\left(1+\rho\right)\left(1+\sigma\right)}C_\alpha\label{SCreqdmlsst2;a},\\
\dfrac{\de C_\alpha}{\de t_\alpha}&=\frac{\sigma}{1+\sigma}{C_\alpha}^2-\frac{\eta+\sigma}{\eta\left(1+\sigma\right)}C_\alpha-\frac{\sigma}{1+\sigma}C_\alpha T_\alpha+\frac{\sigma}{\eta\left(1+\sigma\right)}T_\alpha\label{SCreqdmlsst2;b}. 
\end{align}
\end{subequations}
\end{itemize}
So we have the following scaled problem:
\[\label{TCrdmlss_Cauchy}\tag{$TC\alpha_r$}
\begin{minipage}{.8\textwidth}
\centering 
\textit{Given ${\left[S\right]}_0,{\left[C\right]}_0,{\left[P\right]}_0\,\geq 0$ and $\eta>0$, we seek an interval $\mathcal{I}\subseteq\mathbb{R}$ with $0\in\mathcal{I}$, and a function $\left(T_\alpha,C_\alpha\right)\colon\mathcal{I}\rightarrow {\left[0,1\right]}^2$, such that $\left(T_\alpha,C_\alpha\right)$ satisfies both \eqref{SCreqdmlsst1} if $t_\alpha=\dfrac{t}{t_1}$, or \eqref{SCreqdmlsst2} if $t_\alpha=\dfrac{t}{t_2}$ in ${\left(\mathcal{I}\setminus\left\{0\right\}\right)}^\circ$, and $\left(T_\alpha,C_\alpha\right)=\left(\dfrac{{\left[T\right]}_0}{A_1},\dfrac{{\left[C\right]}_0}{A_1}\right)$, for $t_\alpha=0$.}
\end{minipage}
\]
Working as with problem  (\ref{SC2dmlss_Cauchy}),  we conclude, now for problem (\ref{TCrdmlss_Cauchy}), the following: 
\begin{itemize}
\item (\ref{SCreqdmlsst1;a}) gives that $$
\frac{\de T_\alpha}{\de t_\alpha}\approx 0,\text{ when }\exists B>0\text{ independent of }\eta,\text{ such that: }\left|t_\alpha\right|\leq B,$$
and due to the initial condition of (\ref{TCrdmlss_Cauchy}) we have that $$T_\alpha\approx \frac{{\left[T\right]}_0}{A_1},\text{ when }\exists B>0:\text{ }\left|t_\alpha\right|\leq B.$$ If we insert the above approximate equality into (\ref{SCreqdmlsst1;b}), which will now have the approximate linear form $$\dfrac{\de C_\alpha}{\de t_\alpha}\approx-\frac{\sigma}{1+\sigma}C_\alpha+\frac{\sigma}{1+\sigma}T_\alpha,\text{ when }\exists B>0\text{ independent of }\eta,\text{ such that: }\left|t_\alpha\right|\leq B,$$ then we will get $$C_\alpha\approx T_\alpha+\left(\frac{{\left[C\right]}_0}{A_1}-T_\alpha\right)\exp{\left\{-\frac{\sigma}{1+\sigma}\, t_\alpha\right\}},\text{ when }\exists B>0:\text{ }\left|t_\alpha\right|\leq B,$$ given the initial condition (\ref{TCrdmlss_Cauchy}).\\ Therefore, the initial condition, $\left({\left[T\right]}_{{\rm in}},{\left[C\right]}_{{\rm in}}\right)$, of $\left(\left[T\right],\left[C\right]\right)$ is  $$\left({\left[T\right]}_{{\rm in}},{\left[C\right]}_{{\rm in}}\right)=\left({\left[T\right]}_0,{\left[T\right]}_0+\left({\left[C\right]}_0-{\left[T\right]}_0\right)\exp{\left\{-k_1 A_2 t\right\}}\right).$$
\item From (\ref{SCreqdmlsst2;b}) we obtain that
\begin{equation*}
\begin{split}
\eta\left(\frac{\de C_\alpha}{\de t_\alpha}-\frac{\sigma}{1+\sigma}{C_\alpha}^2+\frac{1}{1+\sigma}C_\alpha-\frac{\sigma}{1+\sigma}C_\alpha T_\alpha\right)=\frac{\sigma}{1+\sigma}\left(-C_\alpha+T_\alpha\right)\approx 0,&\\ &\hspace{-7cm}\text{when }\exists B>0\text{ independent of }\eta,\text{ such that: }\left|t_\alpha\right|\leq B,
\end{split}
\end{equation*}
i.e., $$T_\alpha\approx C_\alpha,\text{ when }\exists B>0\text{ independent of }\eta,\text{ such that: }\left|t_\alpha\right|\leq B.$$ If we insert the above approximate equality in (\ref{SCreqdmlsst2;a}), then the later becomes an approximate linear differential  equation, which is none other than $$\frac{\de T_\alpha}{\de t_\alpha}\approx -\frac{1}{\left(1+\rho\right)\left(1+\sigma\right)}T_\alpha,\text { when }\exists B>0\text{ independent of }\eta,\text{ such that: }\left|t_\alpha\right|\leq B,$$ the solution of which is $$T_\alpha\approx \ell
\exp{\left\{-\frac{1}{(1+\rho)(1+\sigma)}\, t_\alpha\right\}},\text{ when }\exists B>0:\text{ }\left|t_\alpha\right|\leq B,$$ where $\ell\geq 0$ a constant that remains to be determined. Therefore, the outer approximation, $\left({\left[T\right]}_{{\rm out}},{\left[C\right]}_{{\rm out}}\right)$, of $\left(\left[T\right],\left[C\right]\right)$ is 
\begin{equation*}
\begin{split}
\left({\left[T\right]}_{{\rm out}},{\left[C\right]}_{{\rm out}}\right)=\left(\ell A_1\exp{\left\{-k_2 t\right\}},\ell A_1\exp{\left\{-k_2 t\right\}}\right).
\end{split}
\end{equation*}
\end{itemize}

Finally, as usually, we find that $$\ell=\frac{{\left[T\right]}_0}{A_1},$$ as well as that the uniform approximation, $\left({\left[T\right]}_{{\rm un}},{\left[C\right]}_{{\rm un}}\right)$, of $\left(\left[T\right],\left[C\right]\right)$ is 
\begin{equation}
\label{TuCur}
\left({\left[T\right]}_{{\rm un}},{\left[C\right]}_{{\rm un}}\right)=\left({\left[T\right]}_0\exp{\left\{-k_2 t\right\}},{\left[T\right]}_0\exp{\left\{-k_2 t\right\}}+\left({\left[C\right]}_0-{\left[T\right]}_0\right)\exp{\left\{-k_1 A_2 t\right\}}\right).
\end{equation}

\subsection{Conclusions} 

Although for the previous analysis it was used that  $0<\eta\approx 0,$ i.e., $$A_1\ll A_2,$$ nevertheless we emphasise that also the relation $$K_{\!M}\ll A_2,$$ although somehow \textquote{obscure}, plays an essential role in distinguishing (\ref{rQSSA}) from (\ref{sQSSA}). Indeed, let $$K_{\!M}\nll A_2,$$ i.e., $$K_{\!M}\gg A_2\text{ or }K_{\!M}\approx A_2.$$ Then, given $A_1\ll A_2$, we could finally have that $$A_2\ll K_{\!M}+A_1,$$ i.e., (\ref{sQSSA}), if in addition it holds that $K_{\!M}\gg A_2\gg A_1$, or otherwise that $$A_2\approx K_{\!M}+A_1,$$ which does not fall under any case, if additionally $K_{\!M}\approx A_2\gg A_1$.

In addition, we showed that given (\ref{rQSSA}) there is $t_1>0$ such that $$\frac{\de \left[S\right]}{\de t}\approx -k_1 A_2 \left[S\right],\text{ when }\exists B>0\text{ independent of }\eta,\text{ such that: }\left|t\right|\leq B t_1,$$ which arises directly from (\ref{deSalpha2t1approx}), as well as that there is $t_2\gg t_1$ such that $$\frac{\de \left[S\right]}{\de t}\approx \left[S\right] \approx 0,\text{ when }\exists B>0\text{ independent of }\eta,\text{ such that: }\left|t\right|\leq B t_2,$$ which in turn results from (\ref{deSalphat2approxeta}). i.e., we can conclude that $$\upsilon\approx\begin{cases}
k_1 A_2 \left[S\right],&\text{ when }\exists B>0:\text{ }\left|t\right|\leq B t_1,\\
0,&\text{ when }\exists B>0:\text{ }\left|t\right|\leq B t_2,
\end{cases}$$ where $\upsilon$ stands for the rate of the chemical reaction with chemical equation (\ref{SEP}), as we have already mentioned.

\begin{figure}[!htbp]
\centering
\includegraphics[width=.5\textwidth]{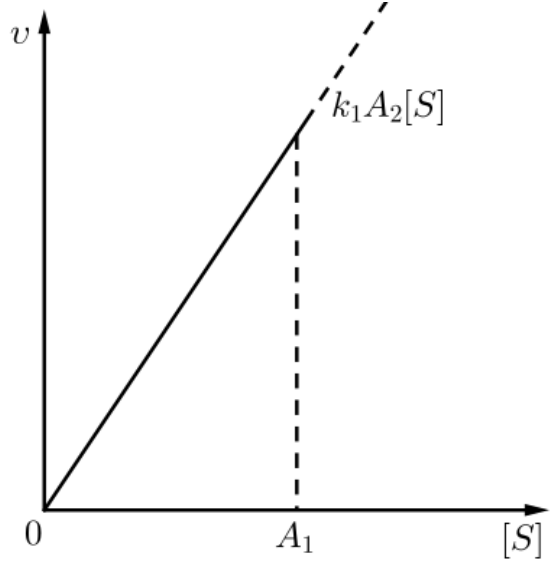}
\caption{An approximation for the kinetics of the chemical reaction (\ref{SEP}) given that (\ref{rQSSA}) holds, for times comparable to $t_1$.}
\label{rQSSA_line}
\end{figure}

Finally, comparing the approximate solution of the free substrate $\left({\left[S\right]}_{{\rm un}},{\left[C\right]}_{{\rm un}}\right)$ of (\ref{SuCur}) with the approximate solution of the total substrate $\left({\left[T\right]}_{{\rm un}},{\left[C\right]}_{{\rm un}}\right)$ of (\ref{TuCur}), we easily observe by definition in (\ref{TclneqSC})  that they are in agreement. 

\subsection{Numerical solution} 

\begin{figure}[!htbp]
\begin{subfigure}{.5\textwidth}
  \centering
\includegraphics[width=1\linewidth]{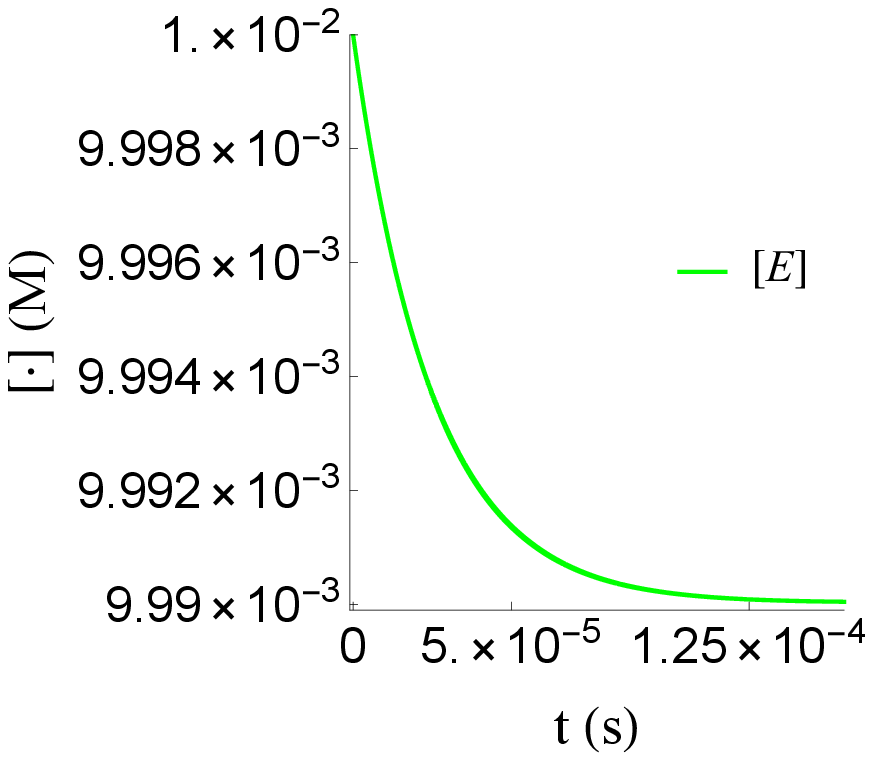}
 \caption{}
  \label{rQSSA-numer-sol.a}
\end{subfigure}
\begin{subfigure}{.5\textwidth}
  \centering
 \includegraphics[width=1\linewidth]{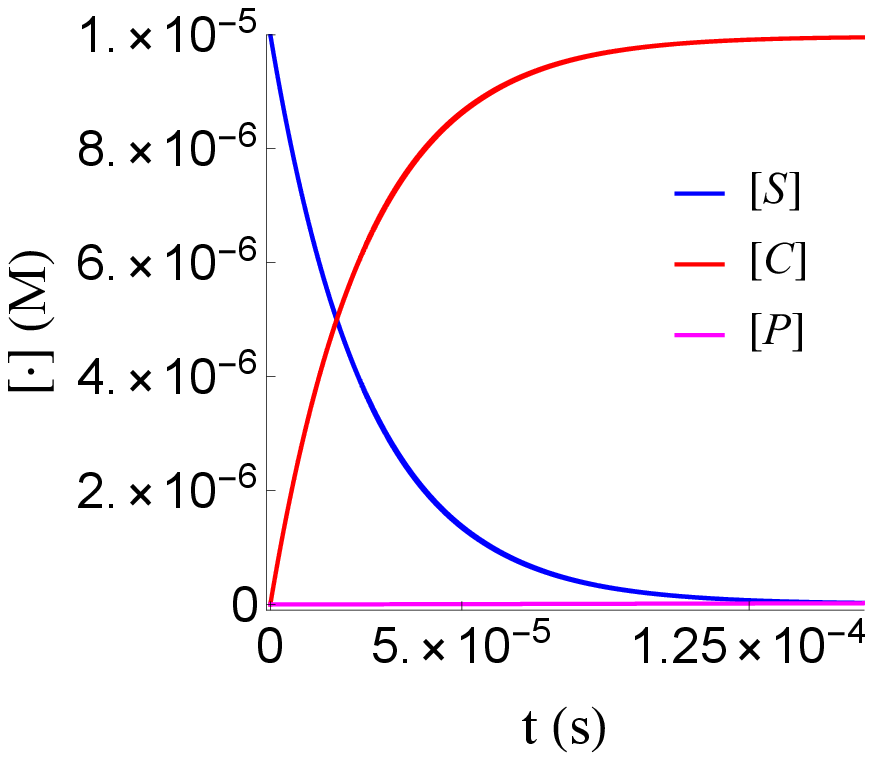}
  \caption{}
  \label{rQSSA-numer-sol.b}
\end{subfigure}
\begin{subfigure}{.5\textwidth}
  \centering
   \includegraphics[width=1\linewidth]{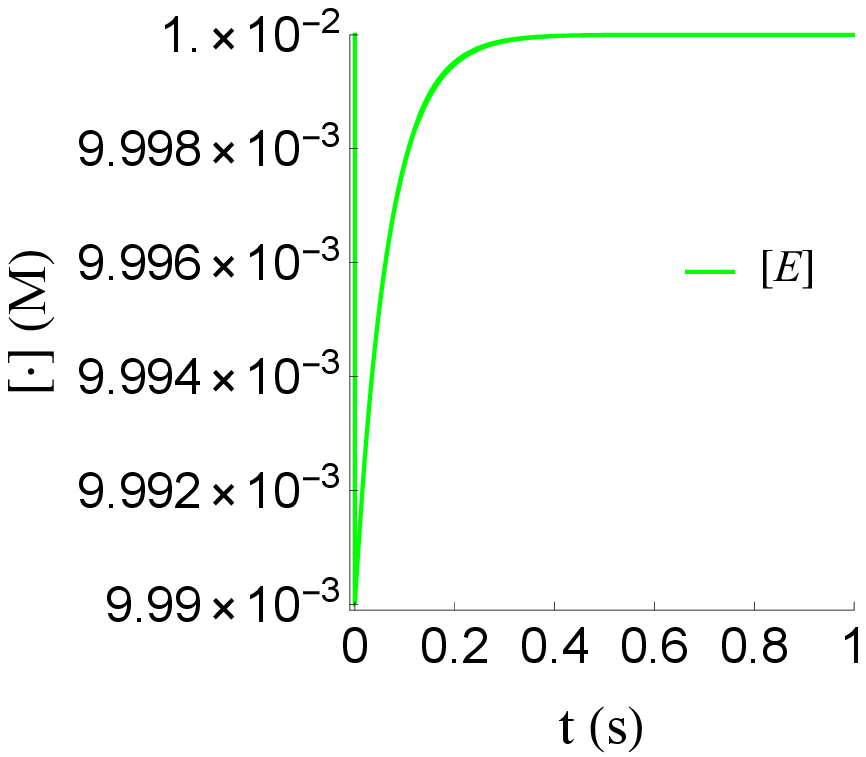}  
  \caption{}
  \label{rQSSA-numer-sol.c}
\end{subfigure}
\begin{subfigure}{.5\textwidth}
  \centering
\includegraphics[width=1\linewidth]{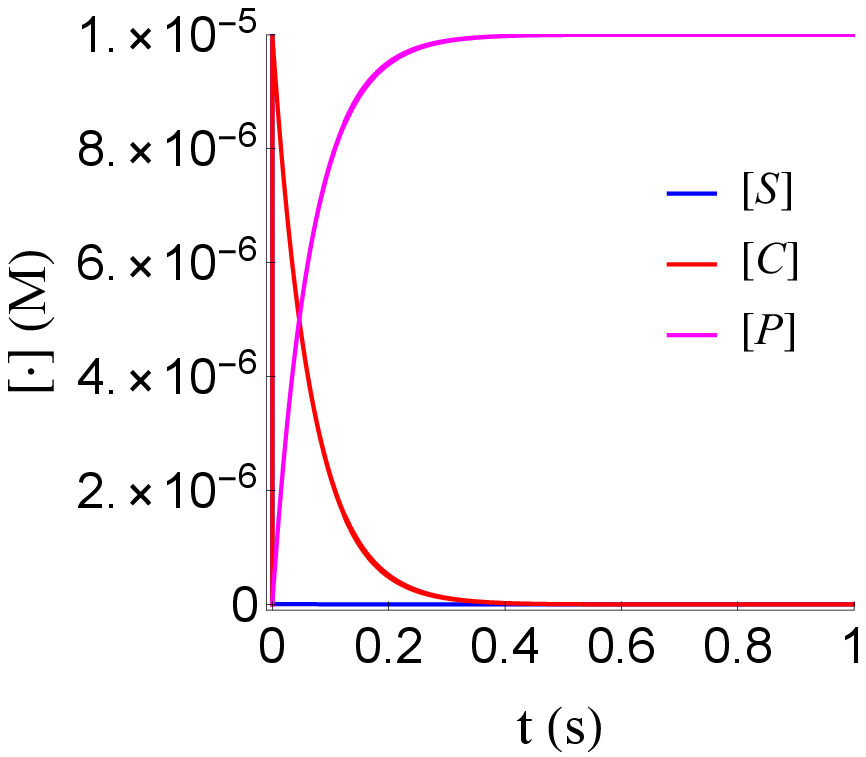}  
  \caption{}
  \label{rQSSA-numer-sol.d}
\end{subfigure}
  \caption{Plots of $\left[S\right]$, $\left[E\right]$, $\left[C\right]$ and $\left[P\right]$ of problem (\ref{SECP_Cauchy}) for non negative times, given that \eqref{rQSSA} holds. We see that $\left[S\right]$ and $\left[C\right]$ are of the same order of magnitude, as well as that there are two distinct phases of the evolution of the phenomenon.}
   \label{rQSSA-numer-sol}
\end{figure}

\begin{figure}[!htbp]
\begin{subfigure}{.5\textwidth}
  \centering
\includegraphics[width=1\linewidth]{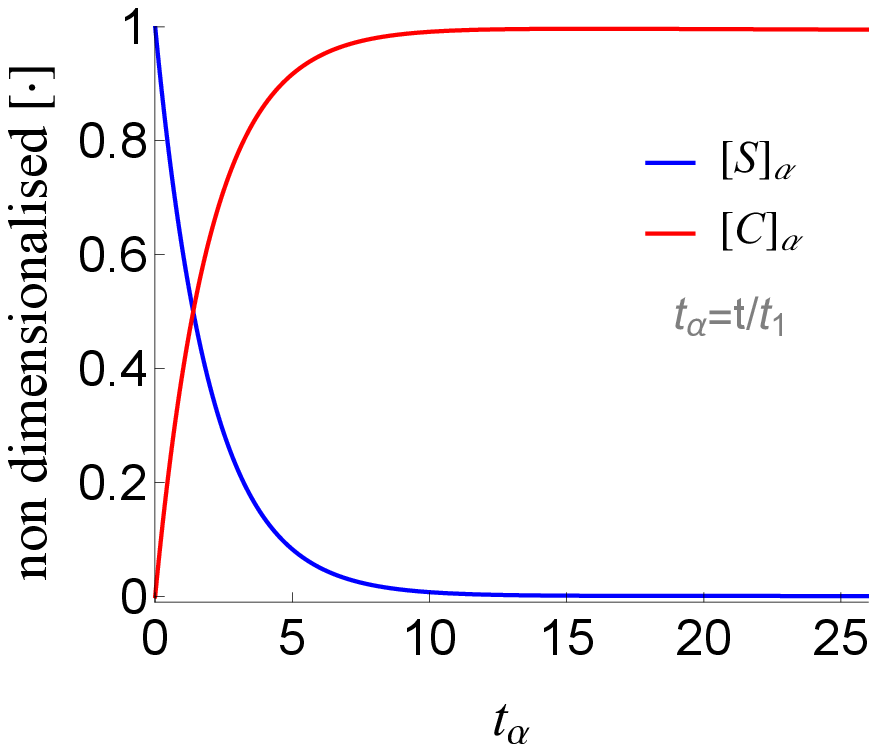}
 \caption{}
  \label{rQSSA-numer-sol-norm.a}
\end{subfigure}
\begin{subfigure}{.5\textwidth}
  \centering
 \includegraphics[width=1\linewidth]{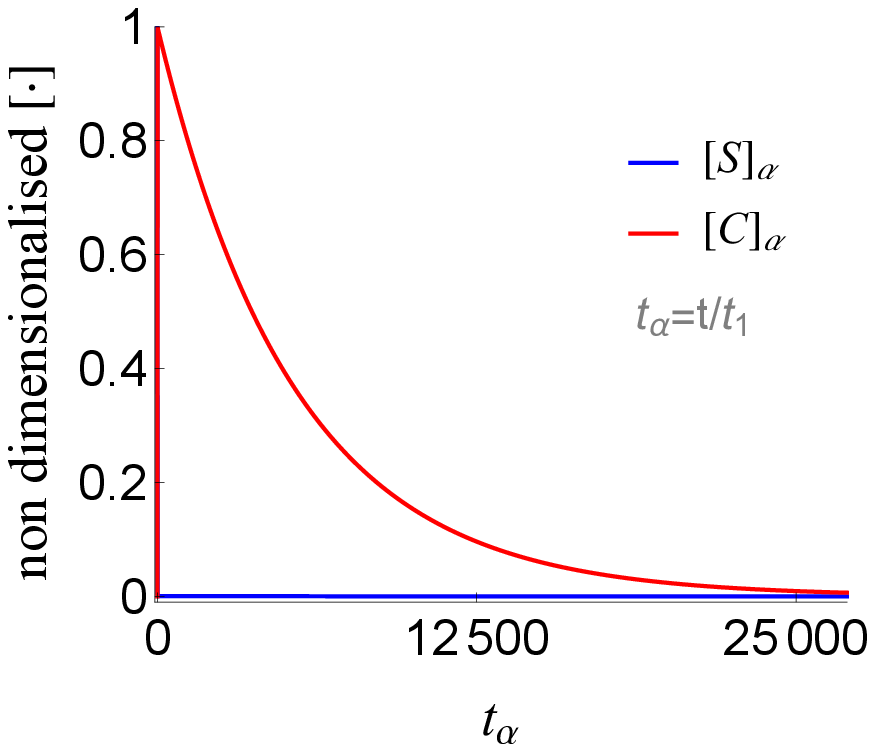}
  \caption{}
  \label{rQSSA-numer-sol-norm.b}
\end{subfigure}
\begin{subfigure}{.5\textwidth}
  \centering
   \includegraphics[width=1\linewidth]{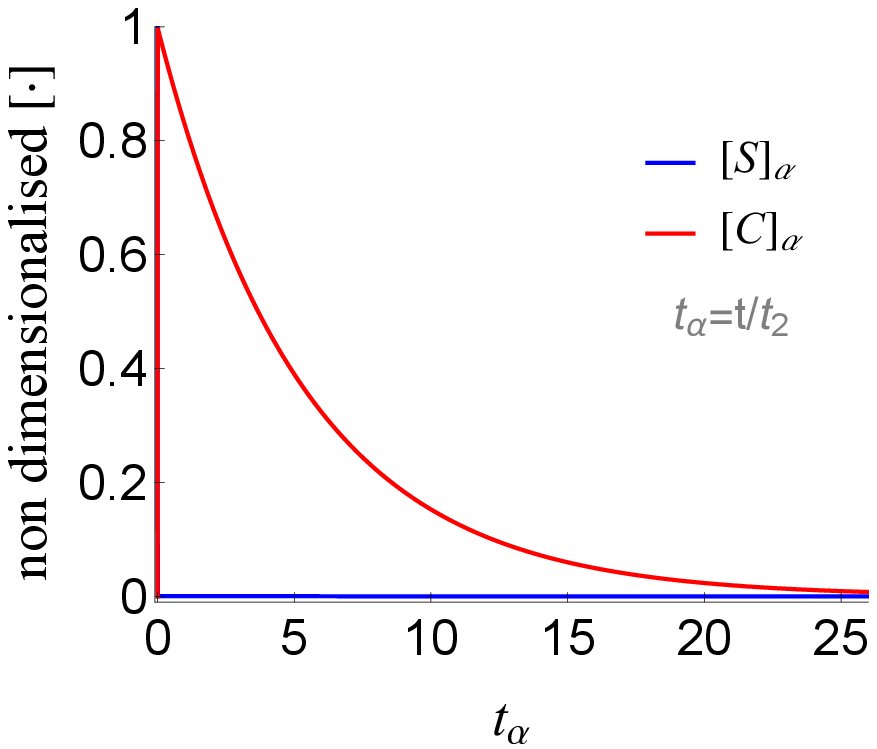}  \caption{}
  \label{rQSSA-numer-sol-norm.c}
\end{subfigure}
\begin{subfigure}{.5\textwidth}
  \centering
\includegraphics[width=1\linewidth]{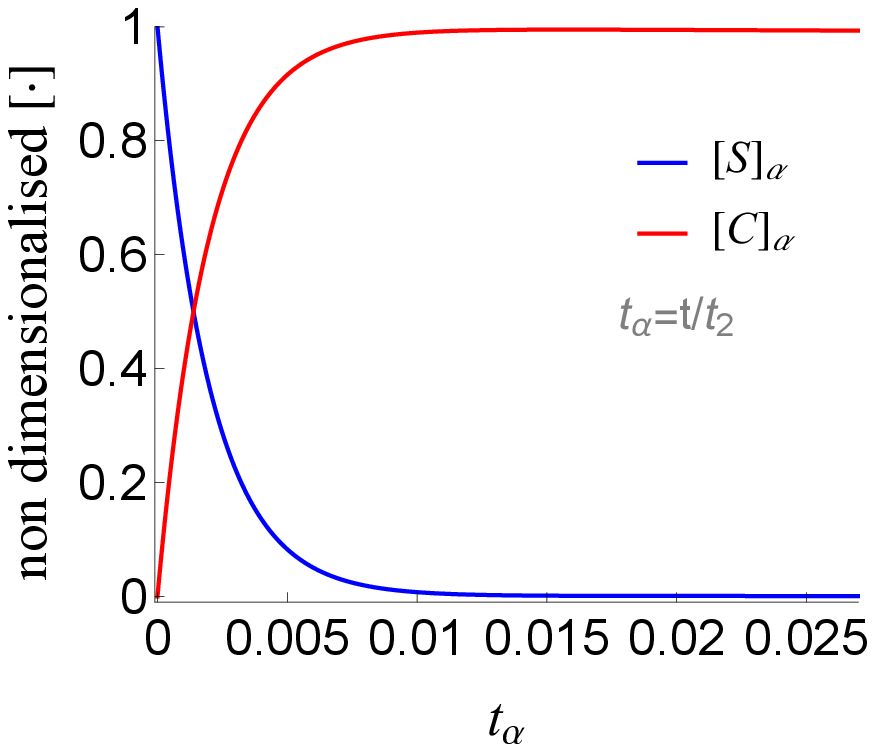}   \caption{}
  \label{rQSSA-numer-sol-norm.d}
\end{subfigure}
  \caption{Plots of $S_\alpha$ and $C_\alpha$ of problem (\ref{SC2dmlss_Cauchy}) for non negative times, given that  \eqref{rQSSA} holds. In (a) and (b) time is measured based on $t_1$, whereas (c) and (d) based on $t_2$.}
   \label{rQSSA-numer-sol-norm}
\end{figure}

\begin{figure}[!htbp]
\begin{subfigure}{.5\textwidth}
  \centering
\includegraphics[width=1\linewidth]{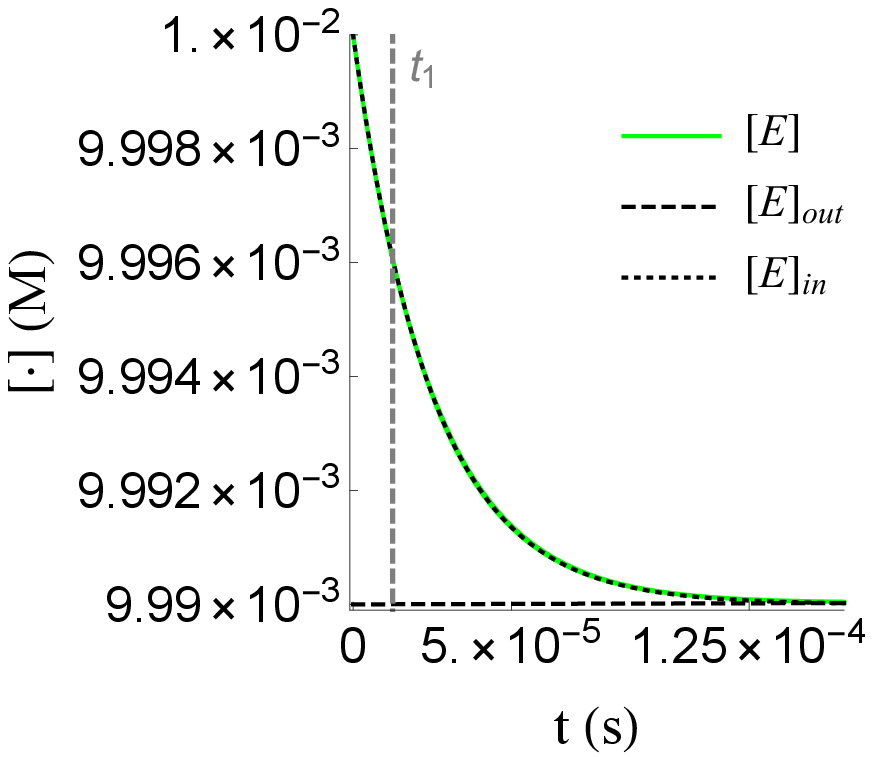}
 \caption{}
  \label{rQSSA-numer-sol-approx.a}
\end{subfigure}
\begin{subfigure}{.5\textwidth}
  \centering
 \includegraphics[width=1\linewidth]{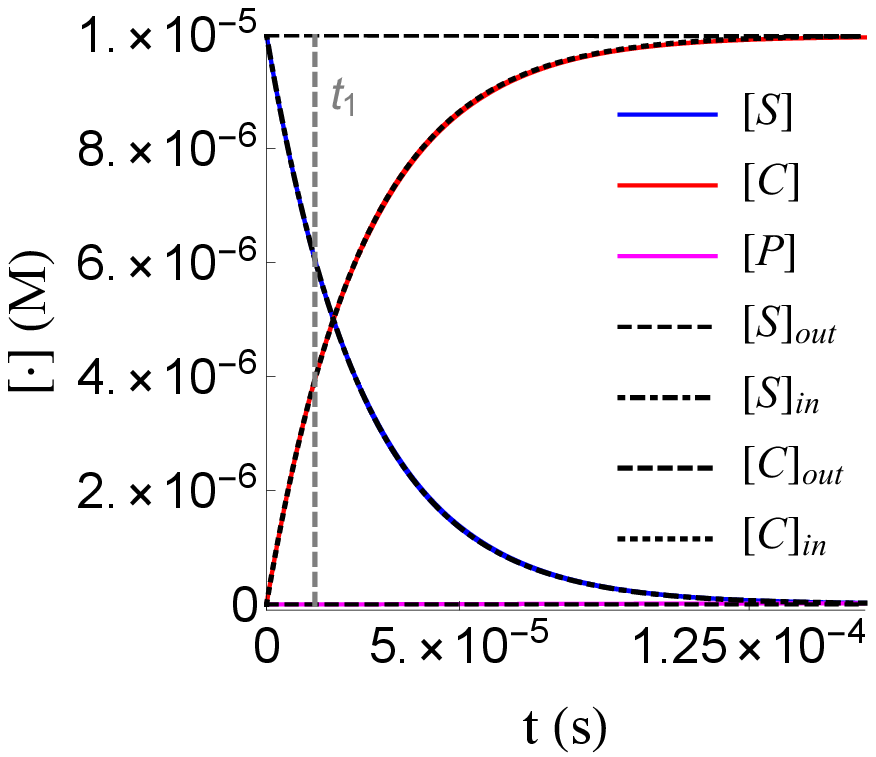}
  \caption{}
  \label{rQSSA-numer-sol-approx.b}
\end{subfigure}
\begin{subfigure}{.5\textwidth}
  \centering
   \includegraphics[width=1\linewidth]{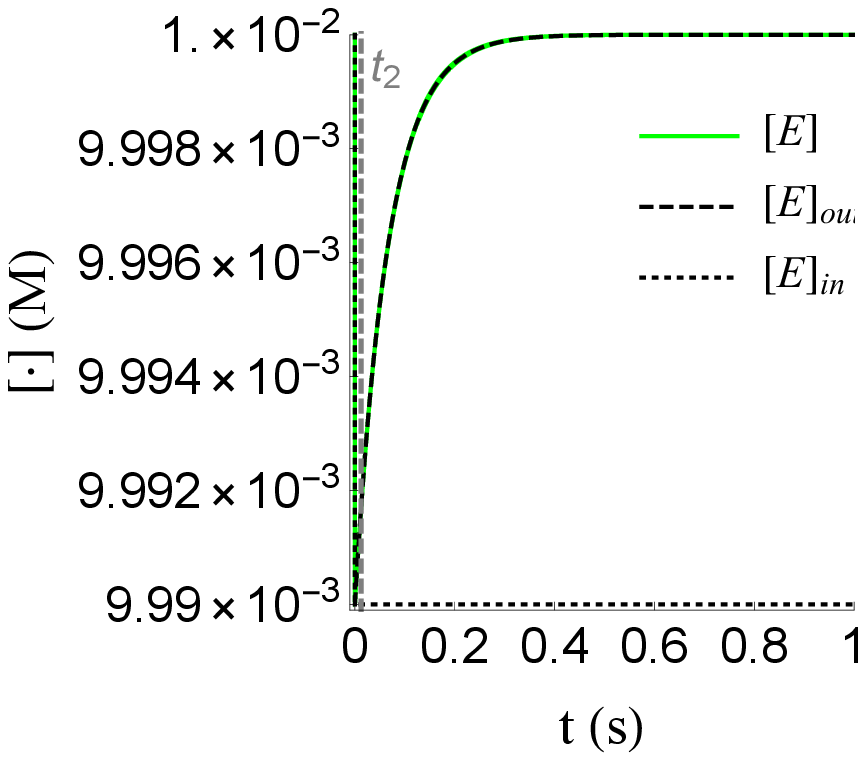} 
  \caption{}
  \label{rQSSA-numer-sol-approx.c}
\end{subfigure}
\begin{subfigure}{.5\textwidth}
  \centering
\includegraphics[width=1\linewidth]{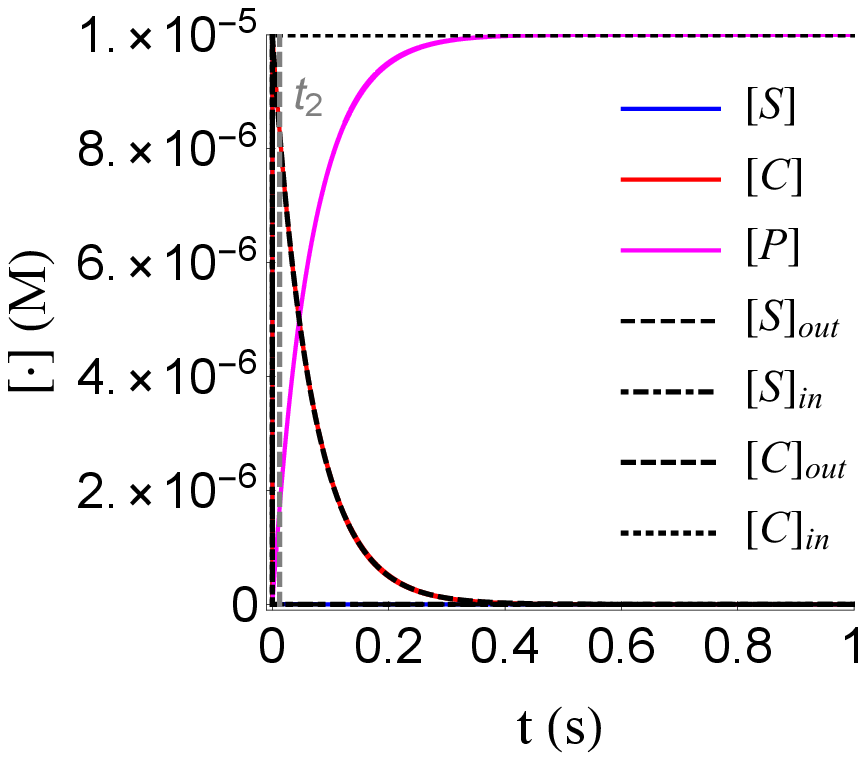}  
  \caption{}
  \label{rQSSA-numer-sol-approx.d}
\end{subfigure}
  \caption{Plots of the inner and outer approximations of ${\left[S\right]}$, $\left[E\right]$ and $\left[C\right]$ of problem (\ref{SECP_Cauchy}), for non negative times, given that \eqref{sQSSA} holds. The inner and outer approximation of  ${\left[E\right]}$ are given by the relations ${\left[E\right]}_{{\rm in}}=A_2-{\left[C\right]}_{{\rm in}}$ and ${\left[E\right]}_{{\rm out}}=A_2-{\left[C\right]}_{{\rm out}}$, respectively, due to (\ref{EplsC}).}
   \label{rQSSA-numer-sol-approx}
\end{figure}

We proceed to the numerical solution of the problem, as shown in \hyperref[rQSSA-numer-sol]{Figure \ref*{rQSSA-numer-sol}}, in \hyperref[rQSSA-numer-sol-norm]{Figure \ref*{rQSSA-numer-sol-norm}} and in \hyperref[rQSSA-numer-sol-approx]{Figure \ref*{rQSSA-numer-sol-approx}}, to verify our conclusions. The numerical values of the constants and the initial conditions are given in the table below.

\begin{table}[ht]
\centering 
\begin{tabular}{c c c}
\hline
Parameter & Value & Unit \\ [0.5ex]
\hline\hline          
$k_{-1}$  & $25$ & $s^{-1}$ \\
$k_1$  & $4\cdot 10^6$ & $M^{-1}s^{-1}$ \\
$k_2$  & $15$ & $s^{-1}$ \\
${\left[S\right]}_0$  & $10^{-5}$ & $M$\\
${\left[E\right]}_0$  & $10^{-2}$ & $M$ \\
${\left[C\right]}_0$  & $0$ & $M$ \\
${\left[P\right]}_0$  & $0$ & $M$ \\  [1ex]      
\hline
\end{tabular}
\label{table:3}
\end{table}
\noindent We calculate $$K_{\!M} ={\left[S\right]}_0=A_1\text{ and }A_2={\left[E\right]}_0=10^3{\left[S\right]}_0=10^3A_1,$$ i.e., $$\varepsilon=500,\text{ }\frac{1}{\varepsilon}=2\cdot 10^{-3}\text{ and }\eta=10^{-3}$$ and $$t_1=1,25\cdot 10^{-5}s\text{ and }t_2=1,25\cdot 10^{-2}s.$$
\section{Discussion} 
\label{Discussion}

We employed a simple scaling algorithm for the rigorous treatment of the quasi-steady-state assumption.

We note that such an algorithm can be utilised in every problem with non negative solutions in a bounded domain, e.g., the simple  classical $SIR$ problem of Epidemiology   
\begin{align*}
\dfrac{\de S}{\de t}&=-\beta SI,\\
\dfrac{\de I}{\de t}&=-\gamma I+\beta SI,\\
\dfrac{\de R}{\de t}&=\gamma I,
\end{align*}
for $\beta,\gamma\,>0$, where the feasible region is the set $$\left\{\left(s,i,r\right)\in{\left[0,N_0\right]}^3\,\big|\,s+i+r=N_0\right\},$$ hence every dependent variable, $S$, $I$ and $R$, is scaled by $N_0$ and the above system then becomes 
\begin{align*}
\frac{1}{\gamma}\dfrac{\de S_\alpha}{\de t}&=-\dfrac{\beta N_0}{\gamma}S_\alpha I_\alpha,\\
\frac{1}{\gamma}\dfrac{\de I_\alpha}{\de t}&=-I_\alpha+\dfrac{\beta N_0}{\gamma}S_\alpha I_\alpha,\\
\frac{1}{\gamma}\dfrac{\de R_\alpha}{\de t}&=I_\alpha.
\end{align*}
By such an approach we naturally obtain the time scale to be $\dfrac{1}{\gamma}$ and, using the well known non dimensionalised quantity $\mathcal{R}_0=\dfrac{\beta N_0}{\gamma}$, the fully scaled equations finally get the form 
\begin{align*}
\dfrac{\de S_\alpha}{\de t_\alpha}&=-\mathcal{R}_0 S_\alpha I_\alpha,\\
\dfrac{\de I_\alpha}{\de t_\alpha}&=-I_\alpha+\mathcal{R}_0S_\alpha I_\alpha,\\
\dfrac{\de R_\alpha}{\de t_\alpha}&=I_\alpha.
\end{align*}
\vspace{0.1cm}

Returning to our problem, from the basic mathematical analysis of (\ref{SECP_Cauchy}) we were able to sleekly generate the quantity $\varepsilon$ that characterises both (\ref{sQSSA}) and (\ref{rQSSA}). Moreover, we naturally determined two pairs of distinctive time scales, each pair of which is characteristic for one of the aforementioned two assumptions.

\begin{figure}[!htbp]
\begin{subfigure}{.5\textwidth}
  \centering
\includegraphics[width=1\linewidth]{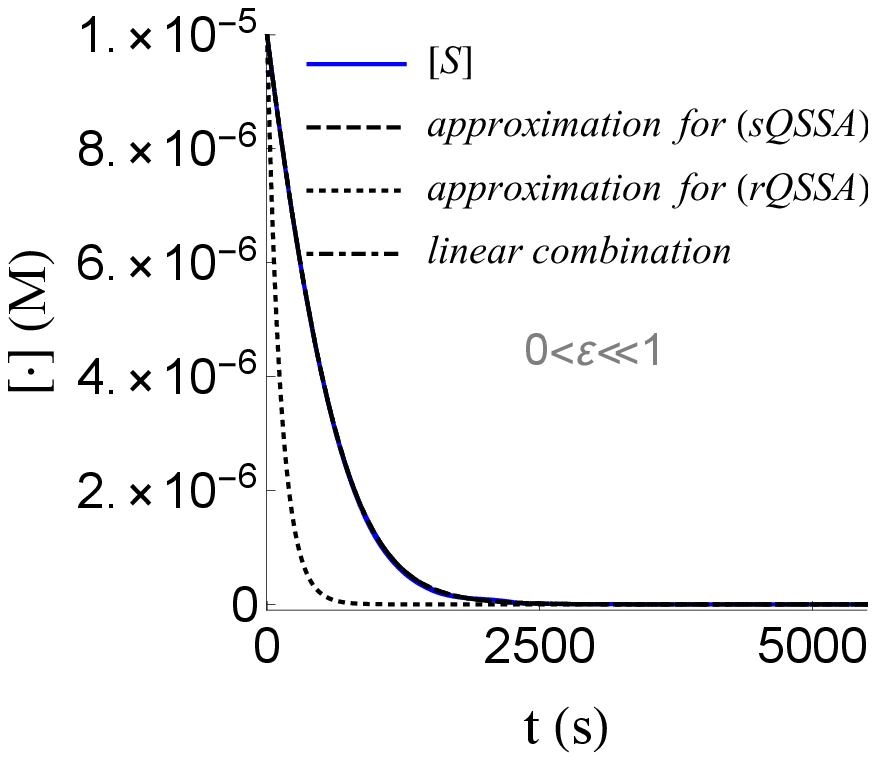}
 \caption{}
  \label{Sgen.a}
\end{subfigure}
\begin{subfigure}{.5\textwidth}
  \centering
 \includegraphics[width=1\linewidth]{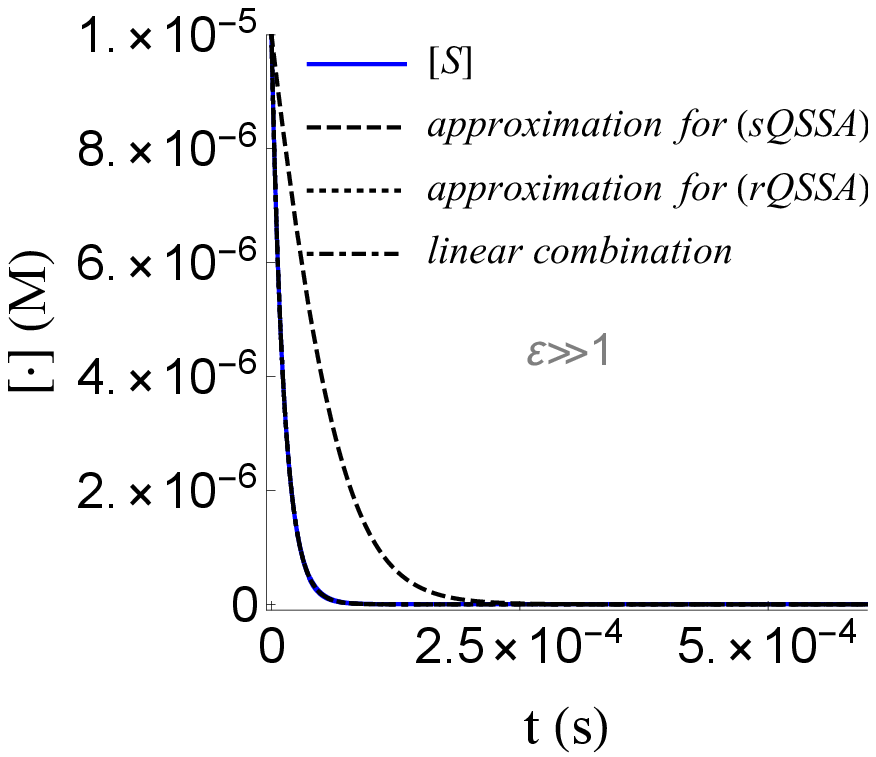}
  \caption{}
  \label{Sgen.b}
\end{subfigure}
\begin{center}
\begin{subfigure}{.5\textwidth}
  \centering
   \includegraphics[width=1\linewidth]{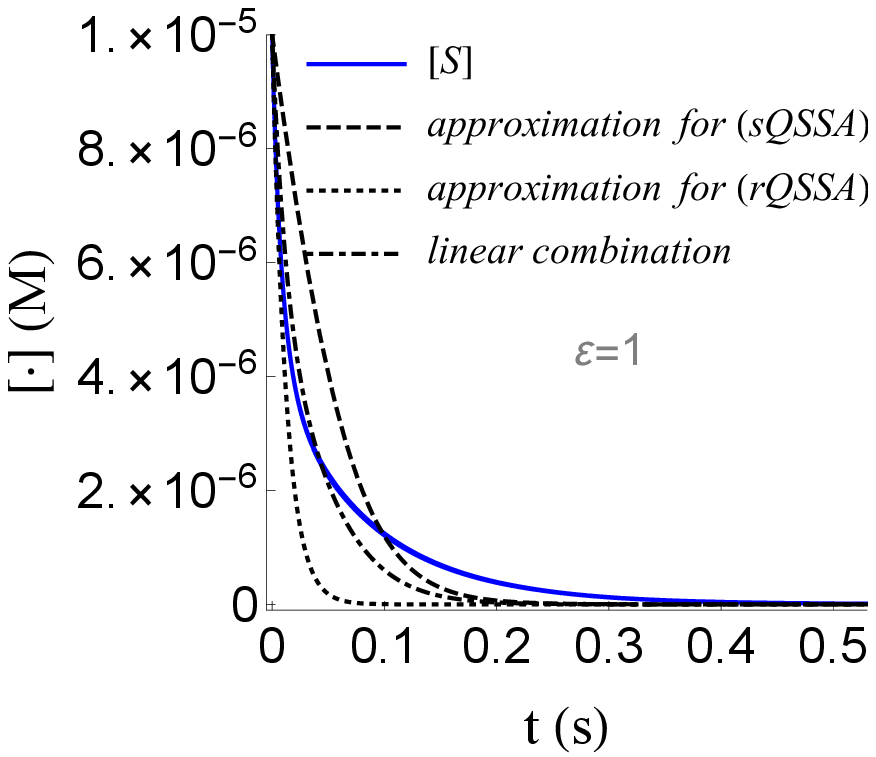} 
  \caption{}
  \label{Sgen.c}
\end{subfigure}
\end{center}
\caption{For a good approximation of the solution of (\ref{SECP_Cauchy}) for the case where $\varepsilon\sim 1$, a sophisticated extrapolation technique is required (work by the present authors in progress), than just a linear combination of the approximations of the solution for (\ref{sQSSA}) and (\ref{rQSSA}).}
   \label{Sgen}
\end{figure}

We further obtained a good approximation of the solution in closed form, for both the cases where $\varepsilon\to 0^+$ and $\varepsilon\to\infty$, which we can communally write as $$\frac{\varepsilon}{1+\varepsilon} \times \left(\text{approximation for }\varepsilon\to\infty\right)\,+\,\left(1-\frac{\varepsilon} {1+\varepsilon}\right) \times \left(\text{approximation for }\varepsilon\to 0^+\right).$$ We emphasise that the above linear combination is far from being a good approximation of the solution for the case where $\varepsilon\sim 1$, as it is illustrated in \hyperref[Sgen]{Figure \ref*{Sgen}}. Such an approximation requires a much more sophisticated extrapolation technique,  the study of which lies beyond the scope of the present work. 
 \bibliographystyle{plain}
\bibliography{mybibfile}\label{bibliography}
\addcontentsline{toc}{chapter}{Bibliography}

\end{document}